\documentclass[a4paper,openbib,twoside,notitlepage]{article}  

\usepackage{amssymb}
\usepackage{latexsym,  pictex}
\usepackage{exscale}
\usepackage{epsf}
\usepackage{theorem}
\usepackage[dvips]{graphicx}

\newfont{\bsym}{cmbsy10 scaled\magstep2}
\newfont{\bsymi}{cmbsy10}
\newfont{\bmath}{cmmib10 scaled\magstep2}
\newfont{\titfont}{cmbx10 scaled \magstep3}
\newfont{\eighrm}{cmr8}
\newfont{\ack}{cmr10}
\newfont{\Ack}{cmbx10 at 14pt}

\newtheorem{theorem}{Theorem}[section]
\newtheorem{lemma}{Lemma}[section]
\newtheorem{prop}{Proposition}[section]
\newtheorem{definition}{Definition}[section]
{\theorembodyfont{\rmfamily} 
\newtheorem{remark}{Remark}[section]}

\newcommand{\vsc}{\vspace{6pt}}
\newcommand{\vsp}{\vspace{12pt}}
\newcommand{\vs}{\vspace{24pt}}

\newcommand{\R}{{\mathbb{R}}}
\newcommand{\NN}{\mathbb{N}}
\newcommand{\LL}{{\bf L}}
\newcommand{\ve}{\varepsilon}
\newcommand{\vfi}{\varphi}

\newcommand{\U}{{\cal U}}

\newcommand{\tv}{\hbox{Tot.Var.}}
\newcommand{\N}{{\cal N}}

\newcommand{\C}{{\cal C}}

\newcommand{\eps}{\varepsilon}

\newcommand{\D}{{\cal D}}
\newcommand{\Z}{\mathbb{Z}}

\newcommand{\n}{\noindent}

\newcommand{\be}{\begin{equation}}
\newcommand{\ee}{\end{equation}}

\newcommand{\lu}{{{\bf L}^1}}

\newcommand{\li}{{{\bf L}^{\infty}}}
\newcommand{\ilu}{{\strut {\bf L}^1}}

\newcommand{\ov}{\overline}

\newcommand{\V}{{\cal V}}

\newcommand{\A}{{\cal A}}

\newcommand{\p}{{\bf p}}

\newcommand{\ess}{{\textrm{ess}\,}}

\begin{document}

\setlength{\voffset}{-0.8in}

\title{
{\huge On the Attainable set for Temple Class Systems}\\
{\huge with Boundary Controls}\\
      \vspace{0.1in}}

\author{
        {\scshape Fabio Ancona}  \thanks{Dipartimento di Matematica and 
                                 C.I.R.A.M., \
                                 P.zza \ Porta \ S. Donato, \ n.~5, \
                                 40123 - Bologna, \ Italy;
                                 \newline 
                                 E-mail: \texttt{ancona@ciram3.ing.unibo.it.
                                 \vspace{6pt}}}
\and    
        {\scshape Giuseppe Maria Coclite} \thanks{SISSA-ISAS,
                                via Beirut 2-4, 34014 - Trieste, \ Italy;      
                                      E-mail: \texttt{coclite@sissa.it.}} 
                                \vspace{15pt}}

\date{May 2002}

\maketitle
\thispagestyle{empty}

%
\begin{abstract}
\vspace{5pt}

\noindent

Consider the initial-boundary value problem for
a strictly hyperbolic, genuinely nonlinear, Temple class system
of conservation laws 
$$
u_t+f(u)_x=0, \qquad\ \
u(0,x)=\ov u(x), \qquad
\left\{\!\!\!\!\!\!\!\!
\begin{array}{ll}
&u(t,a)=\widetilde u_a(t),  \\
\noalign{\smallskip}
&u(t,b)=\widetilde u_b(t), 
\end{array}
\right.
\eqno(1)
$$
on the domain $\Omega =\{(t,x)\in\R^2 : t\geq 0,\, a \le x\leq b\}.$
We study the mixed problem (1) 
from the point of view of control theory,
taking the initial data $\overline u$
fixed, and regarding the boundary data $\widetilde u_a, \, \widetilde u_b$
as control functions
that vary in prescribed sets $\U_a,\, \U_b$, of $\li$ boundary controls.
In particular, we consider
the family of configurations
$$
\A(T\,) \doteq \big\{
u(T,\cdot)~; ~~
u \ {\rm \ is \ a \ sol. \ \ to} \ \ (1),
\quad \widetilde u_a\in \U_a, \ \, \widetilde u_b \in \U_b
\big\}
$$
that can be attained by the system at a given time $T>0$,
and we give a
description of the attainable
set  $\A(T)$
in terms of suitable Oleinik-type conditions.
We also establish closure and compactness of the set $\A(T)$ 
in the $\lu$~topology.
\end{abstract}
\vspace{1.5cm}

\begin{center}
Ref. S.I.S.S.A. 10/2002/M (February 2002)
\end{center}

\vspace{1.5cm}

1991\textit{\ Mathematical Subject Classification:} 
35L65,
35B37

\vfill
\pagebreak
%
%
%

\setlength{\voffset}{-0in}  
\setlength{\textheight}{0.9\textheight}

\thispagestyle{empty}
\null
\vfill
\pagebreak

\thispagestyle{plain}
\setcounter{page}{1}
\setcounter{equation}{0}

\section[]{Introduction} 
\label{section31}
\indent

Consider the initial-boundary value problem for
a strictly hyperbolic, genuinely nonlinear, system
of conservation laws in one space dimension
\begin{eqnarray}
& &u_t+f(u)_x=0,
\label{311} \\
& &u(0,x)=\ov u(x), 
\label{312} \\
& &u(t,a)=\widetilde u_a(t), 
\label{313} \\
& &u(t,b)=\widetilde u_b(t), 
\label{314}
\end{eqnarray}
on the strip $\Omega =\{(t,x)\in\R^2 ~;~ t\geq 0,\, x\in[a,\,b]\}$.  
Here, $u=u(t,x)\in \R^n$ is the vector of the
conserved quantities, $\widetilde u_a, \, \widetilde u_b$
are measurable, bounded boundary data, and 
the flux function
$f:U \mapsto\R^n$ is a smooth vector field defined on some open set
$U \subseteq \R^n$, that belongs to a class of fields
introduced by Temple \cite{temple, serre} for which rarefaction and Hugoniot
curves coincide.
We recall that, for problems of this type,
classical solutions may develop discontinuities in finite time,
no matter of the regularity of the initial and boundary data.
Hence, it is natural to consider weak solutions in the sense 
of distributions.
Moreover,  
since, in general, the Dirichlet 
conditions (\ref{313})-(\ref{314}) cannot be fulfilled pointwise a.e.
(see \cite{barlerned, dubfiocco}),
different weaker formulations of the boundary 
condition have been 
considered in the literature
(see \cite{am, jofiocco, sable}  and
references therein). Here, 
following F.~Dubois, P.G.~LeFloch \cite{dubfiocco},
we will adopt a formulation of (\ref{313})-(\ref{314}) 
based on the definition of a time-dependent
set of {\it admissible
boundary data}, that is related to the notion of Riemann problem.
\vsc

In the present paper, having in mind applications of Temple systems to 
problems of oil reservoir simulation, 
multicomponent chromatography, as well as in
traffic flow models, we study the effect of the boundary conditions
(\ref{313})-(\ref{314}) on the solution
of (\ref{311})-(\ref{312}) from the point of view of
control theory. Namely, following the same approach
adopted in \cite{ancmar1, ancmar2} for scalar conservation laws,
we fix an initial data $\overline u\in \li([a, b])$
and we consider
the family of configurations
\be
\A(T;\, \U_a,\, \U_b) \doteq \big\{
u(T,\cdot)~; ~~
u \ {\rm \ is \ a \ sol. \ \ to} \ \ (\ref{311})-(\ref{314}),
\quad \widetilde u_a\in \U_a, \ \, \widetilde u_b \in \U_b
\big\}
\label{attset1}
\ee
that can be attained at a given time $T>0$
by solutions to  (\ref{311})-(\ref{314}),
with boundary data $\widetilde u_a, \,
\widetilde u_b$ that vary in 
prescribed sets $\U_a, \, \U_b\subset\li(\R^+)$
of {\it admissible boundary controls.}
In the case of scalar, convex conservation laws, 
it was proved in \cite{ancmar1},
by using the theory of generalized characteristics~\cite{daf},
that the profiles $w(x)$ which can be 
attained at a fixed 
time $T>0$ are only 
those for which the map $x \mapsto \frac{f'(w(x)}{x}$ 
is non increasing.
Under the assumption that $f'(u) \geq 0$
for all $u$, and for solutions
of the mixed problem (\ref{311})-(\ref{314}) on the region~$\Omega$,
this condition is equivalent to
the Oleinik-type inequalities
\be
D^+w(x)  \leq \displaystyle{\frac{f'(w(x))}{(x-a)\, f''(w(x))}}
\qquad {\rm for \ \ a.e.} \ \ x\in[a,b]\,,
\ee
($D^+w$ denoting the upper Dini derivative of $w$).
For general $n\times n$ systems,  
a complete
characterization of the attainable set  
does not seem 
possible, due to the complexity of repeated wave-front
interactions. However, in the particular case of Temple systems,
 wave interactions 
can only change the speed of wave-fronts, without 
modifying  their amplitudes, due to
the special geometric
features of such systems.
Therefore,  the only restriction to boundary controllability 
is the decay  due
to genuine nonlinearity.
We then consider here a convex, compact set $\Gamma\subset U$,
and provide a description of the
attainable set 
$$
\A(T)\doteq \A(T;\, \U^\infty,\, \U^\infty),
\qquad\quad \U^\infty\doteq \li([0,T], \Gamma)\,,
$$
in terms of certain Oleinik-type conditions.
We also establish the compactness of 
$\A(T)$
in the $\lu$~topology.
Applications to calculus of variations and problems
of optimization (where the cost functional depends on the
profile of the solution at a fixed time $T$)
motivate the study of topological properties
of $\A(T)$.
\vsp

The paper is organized as follows.
Section~\ref{section32} contains the basic definitions and the
statement of the main results.
We also provide in this section
a review of the existence and well-posedness theory
for the mixed problem~(\ref{311})-(\ref{314}),
and  a description of a front tracking algorithm
that will be used throughout the paper.
In Section~\ref{section33} we establish some preliminary 
estimates, and a regularity result concerning 
the global structure of 
solutions to the mixed problem~(\ref{311})-(\ref{314})
generated by a front tracking algorithm. 
The proof of the main results is contained in Section~\ref{section34}.
\vsp

\section{Preliminaries and statement of the main results} 
\label{section32}
\setcounter{equation}{0}
\vspace{6pt}

\subsection{Formulation of the problem}
\label{subsec:formulation}
\indent

Let $f:U \mapsto\R^n$ be the flux function
of the strictly hyperbolic system (\ref{311})
defined on a neighborhood of the origin $U\subseteq\R^n$.\
Denote by
$\lambda_1(u)<\cdots <\lambda_n(u)$
the eigenvalues of the Jacobian matrix $Df(u)$,
and let $\{r_1(u), \dots , r_n(u)\}$  be a basis of
right eigenvectors
of $Df(u)$.
By possibly considering a 
sufficiently small restriction
of the domain $U$, we may assume
that the following {\it uniform} strict hyperbolicity condition holds.

\begin{itemize}
\itemindent 10pt
\item[{\bf (SH1)}] 
For every $u, v \in U,$ the characteristic speeds at these points
satisfy
\be
\lambda_i(u) < \lambda_j(v),
\qquad\quad
\forall~
1\leq i<j \leq n\,.
\label{unifstrhyp}
\ee
\end{itemize}
We also assume that there is a fixed set
of characteristic lines
entering the interior of the strip $[a,b]\times \R^+$  at
the boundaries~$x=a,\>x=b$, i.e. that,
for some index $p\in\{1,\ldots,n\}$, there holds
\be
\lambda_{p}(u)< 0<\lambda_{p+1}(u),\qquad\quad\forall~u\in U,
\label{321}
\ee
and we let $\lambda^{\min}, \, \lambda^{\max}$ 
denote the minimum and maximum characteristic speed
so that there holds
\be
\label{uppbound}
0<\lambda^{\rm min} \leq |\lambda_i(u)| \leq \lambda^{\rm max},
\qquad\quad \forall~u\in U.
\ee

Moreover, we assume that  each $i$-th characteristic field $r_i$
is {\it genuinely nonlinear} in the sense
of Lax~\cite{lax},
and that
system (\ref{311}) is of Temple class according with the following.
\begin{definition}  \label{definition31}
A system of conservation laws is of {\it Temple class} if there exists
a system of 
coordinates $w=(w_1,\ldots,w_n)$ consisting of Riemann invariants, and such 
that the level sets \textrm{$\big\{u\in U;~w_i(u)=constant\big\}$} 
are hyperplanes (see \cite{serre}). 
\end{definition}

By possibly performing a translation of coordinates,
it is not restrictive to
assume that the Riemann invariants
are chosen so that 
$\partial_i\, \lambda_i(w)>0$, $i=1,\dots n$,
 for all $w=w(u), \, u\in U$.
Throughout the paper, we will often write
$w_i(t,x)\doteq w_i\big(u(t,x)\big)$ to denote the $i$-th
Riemann coordinate of a solution $u=u(t,x)$ to (\ref{311}).
We recall that,
for a Temple class system, 
Hugoniot curves and rarefaction curves coincide~\cite{temple}.
Moreover, as observed in~\cite{A-G},
thanks to the existence of Riemann coordinates
one can show that the assumption {\bf SH1}
implies the invertibility of the map $f : U \mapsto f(U).$
%
%
\vsp

We next introduce a definition of weak solution 
to (\ref{311})-(\ref{314})
which includes an entropy admissibility condition
of Oleinik type on the decay of positive waves,
so to achieve uniqueness.
The boundary conditions (\ref{313})-(\ref{314}) 
are formulated in terms of the weak trace
of the flux $f(u)$ at the the 
boundaries $x=a,\, x=b$, and are related to the notion
of Riemann problem in the same spirit of~\cite{dubfiocco}. 
To this purpose, 
letting $u(t,x)=W(\xi=x/t;\, u_L, u_R),$  \ $u_L,\, u_R\in U$,
denote the self-similar solution of the Riemann problem for (\ref{311}) 
with initial data
$$
u(0,x)=\left\{
\begin{array}{ll}
\!\!u_L  &\textrm{if $x<0$},
\\
\noalign{\smallskip}
\!\!u_R &\textrm{if $x>0$},
\end{array} \right.
$$
for any given boundary state $\widetilde u\in
U$, we define the set of {\it admissible states at the boundaries} 
{\setlength\arraycolsep{2pt}
\be
\begin{array}{ll}
\V_a(\widetilde u)
&\doteq\big\{W(0+;\, \widetilde u, u_R)~;~ u_R\in U\big\},
\\
\noalign{\medskip}
\V_b(\widetilde u)
&\doteq\big\{W(0-;\, u_L, \widetilde u)~;~ u_L\in
U\big\}. 
\end{array}
\label{324} 
\ee } 
\begin{definition} \label{definition32} 
A function $u:[0,\,T]\times [a,b]\mapsto U$  
is an entropy
weak solution of the initial-boundary value problem
(\ref{311})-(\ref{314})
on $\Omega_T\doteq [0,\,T]\times [a,b],$
if it is continuous as a function from $]0,\,T]$ 
into $\lu$,
and the following properties hold:
\begin{itemize}
\item[(i)] $u$ is a distributional solution to the Cauchy problem
(\ref{311})-(\ref{312}) on $\Omega_T$ 
in the \linebreak sense 
that, 
for every  test function $\phi\in\C^1_c$ 
with compact support contained in the set 
$\{(t,x)\in\R^2;\ a<x<b,\, t<T\},$
there holds
$$
\int_0^{T}\int_a^b\big( u(t,x)\cdot\phi_t(t,x)+f(u(t,x))
\cdot\phi_x(t,x)\big)dx~dt+\int_a^b \ov u(x)\cdot\phi(0,x)dx=0\,;
$$ 

\item[(ii)] the flux $f(u)$ admits weak$^*$ traces at the boundaries
$x=a,\>x=b$, i.e. there exist two measurable functions $\Psi_a,\Psi_b  :
[0,\,T] \mapsto\R^n$  
such that
\be
f(u(\cdot,x)) \ \mathop{-\!\!\!-\!\!\!\rightharpoonup}^*_{~x \to a^+} \
\Psi_a,\qquad
f(u(\cdot,x)) \ \mathop{-\!\!\!-\!\!\!\rightharpoonup}^*_{~x \to b^-} \
\Psi_b
\qquad\hbox{in }~~\li([0,\,T]),
\label{328}
\ee
and the boundary conditions (\ref{313})-(\ref{314}) 
are satisfied in the following sense
\be
\Psi_a(t) \in f\big(\V_a(\widetilde u_a(t))\big),
\qquad\Psi_b(t) \in f\big(\V_b(\widetilde u_b(t))\big)
\qquad\hbox{ for ~a.e. \ \ $0\leq t \leq T$}; \label{329}
\ee
\item[(iii)] $u$ satisfies the following entropy conditions
on the decay of positive waves in time and in space.
There exists some constant $C>0,$ depending only on the system (\ref{311}),
so that
\begin{itemize}
\item[(a)]
For any  $0<t\leq T,$ and for a.e. $a<x<y<b,$ there holds
{\setlength\arraycolsep{0pt}
\begin{eqnarray}
w_i(t,y)-w_i(t,x) ~&\leq& ~C\cdot \bigg\{{{y-x}\over
{t}}+\log\bigg(\frac{y-b}{x-b}\bigg)\bigg\}
\quad\ \textrm{if}\quad \ \  i\in\{1,\ldots,p\}, \nonumber
\\
\label{330} 
\\ & {}
&{} \nonumber \\ w_i(t,y)-w_i(t,x) ~&\leq& ~
C\cdot
\bigg\{{{y-x}\over {t}}+\log\bigg(\frac{y-a}{x-a}\bigg)\bigg\} 
\quad\ \textrm{if}\quad \ \ i\in\{p+1,\ldots,n\};
\nonumber \\
\label{331}
\end{eqnarray}}
\item[(b)]
For a.e.  $a<x<b,$ and for a.e. $0<\tau_1<\tau_2\leq T,$
there holds
{\setlength\arraycolsep{0pt}
\begin{eqnarray}
w_i(\tau_2,x)-w_i(\tau_1,x) ~&\leq& ~C\cdot
\bigg\{{{\tau_2-\tau_1}\over {x-b}}+
\log\left(\frac{\tau_2}{\tau_1}\right)\!\!\bigg\}
\quad\ \textrm{if}\quad\ \ i\in\{1,\ldots,p\}, \nonumber
\\
\label{327b} \\
& {} &{} \nonumber \\
w_i(\tau_2,x)-w_i(\tau_1,x) ~&\leq&
~C\cdot\bigg\{{{\tau_2-\tau_1}\over {x-a}}+
\log\left(\frac{\tau_2}{\tau_1}\right)\!\!\bigg\} 
\quad\ \textrm{if}\quad \ \ i\in\{p+1,\ldots,n\}.
\nonumber
\\
\label{328b}
\end{eqnarray}}
\end{itemize}
\end{itemize}
\end{definition}

\begin{remark} \label{remark20a}
The set of admissible flux values at the boundaries $x=a$,
$x=b$,
can be expressed in Riemann coordinates as 
\be
{\setlength\arraycolsep{2pt}
\begin{array}{ll}
f\big(\V_a(\widetilde u)\big)
&=\Big\{f(u)~;~ w_i(u)=w_i(\widetilde u)
 \qquad \forall~i=p+1,\dots, n\Big\},
\\
\noalign{\medskip}
f\big(\V_b(\widetilde u)\big)&=\Big\{f(u)~;~ w_i(u)=w_i(\widetilde u)
 \qquad \forall~i=1,\dots, p\Big\}.
\end{array}}
\label{329c}
\ee
Hence, by the invertibility of the 
map $f: U\mapsto f(U)$,
the above boundary conditions~(\ref{329}) are
equivalent to 
the set of equalities
\be
{\setlength\arraycolsep{1pt}
\begin{array}{ll}
w_i\big(f^{-1}(\Psi_a(t))\big)&=w_i\big(\widetilde u_a(t)\big)
\qquad\hbox{ for ~a.e. \ \ $0\leq t \leq T,$}
\qquad i=p+1,\dots ,n,
\\
\noalign{\medskip}
w_i\big(f^{-1}(\Psi_b(t))\big)&=w_i\big(\widetilde u_b(t)\big)
\qquad\hbox{ for ~a.e. \ \ $0\leq t \leq T,$}
\qquad i=1,\dots ,p.
\end{array}}
\label{3211}
\ee
This means that the boundary conditions (\ref{329}) guarantee that,
at almost every \linebreak time  $t\in [0,\,T],$ 
the solution to the Riemann problem for (\ref{311}), 
having left and right initial states 
$u^L=\widetilde u_a(t),$  $u^R=f^{-1}(\Psi_a(t)),$
contains only waves with negative speeds,
while the solution to the Riemann problem
with initial states
$u^L=f^{-1}(\Psi_b(t)),$  $u^R=\widetilde u_b(t),$
contains only waves with positive speeds.
Thus, in particular, such solutions
do not contain any front entering the domain 
$[t,\,+\infty[\,\times ]a, b[$.
\end{remark}
\vsp

In the present paper we regard the boundary data
as admissible controls and, in connection
with a fixed convex, compact set $\Gamma\subset U$ having the form
\be
\Gamma=\Big\{u\in U; ~~~
w_i(u)\in [\alpha_i,~\beta_i],\qquad i=1,\ldots,n\Big\},
\label{3212}
\ee
we study the basic properties
of the 
{\it attainable set} for (\ref{311})-(\ref{312}),
i.e. of the set
\be
\A(T)\doteq
\big\{
u(T,\cdot)~; ~~
u \ {\rm \ is \ a \ sol. \ \ to} \ \ (\ref{311})-(\ref{314}),
\quad \widetilde u_a,\, \widetilde u_b \in 
{\bf L}^{\infty}([0,T],\Gamma)
\big\}
\label{attset2}
\ee
which consists of all profiles that can be attained 
at a fixed time $T>0,$ by 
entropy weak
solutions of (\ref{311})-(\ref{314})
(according with Definition~\ref{definition32})
with a fixed
initial data $\overline u\in \li([a, b],\Gamma)$, and boundary data 
$\widetilde u_a, \>\widetilde u_b $ that vary in
\be
\U_T^\infty\doteq{\bf L}^{\infty}([0,T],\Gamma).
\label{controlset}
\ee
We will establish a
characterization of (\ref{attset2})
in terms of certain Oleinik type estimates on the decay
of positive waves,
and we will prove the compactness of~(\ref{attset2})
in the $\lu$ topology.
\vspace{5pt}

\subsection{Statements of the main results}
\label{subsec:statement}
\indent

For any $\rho>0$, consider the set of maps 
\be 
\begin{array}{ll}
\!\!\!\!\!&\!\!\!
K^\rho\!\doteq\!
\left\{ \vfi\in{\bf L}^{\infty}([a,b],\,\Gamma)~;\> \,
\begin{array}{ll}
\displaystyle{\frac{w_i(\vfi(y))-w_i(\vfi(x))}{y-x}}
\leq \displaystyle{\rho\over x-a}\quad 
\left\{
\!\!\!\!\!\!\!\!
\begin{array}{ll}
&\hbox{for \, a.e.} \quad a<x<y<b,
\\
\noalign{\smallskip}
&\hbox{if} \quad i\in \{p+1,...,n\}
\end{array}
\right.
\\
\noalign{\medskip}
\displaystyle{\frac{w_i(\vfi(y))-w_i(\vfi(x))}{y-x}}
\leq \displaystyle{\rho\over b-y}\quad 
\left\{
\!\!\!\!\!\!\!\!
\begin{array}{ll}
&\hbox{for \, a.e.} \quad a<x<y<b,
\\
\noalign{\smallskip}
&\hbox{if} \quad i\in \{1,...,p\}
\end{array}
\right.
\end{array}
\!\!\!\!\!
\right\}.
\\
\noalign{\medskip}
\end{array}
\label{K}   
\ee
The inequalities in (\ref{K}) reflect the fact that positive waves
entering through the boundaries $x=a, \, x=b$ decay in time.
Therefore, their density (expressed in terms of Riemann coordinates)
is inversely proportional to their distance
from their entrance point on the boundary.   

\begin{theorem} \label{theorem1}
Let  (\ref{311}) be a system of Temple class
with all characteristic fields genuinely nonlinear, and assume
that the strict hyperbolicity condition ${\bf (SH1)}$ is
verified.
Then, for every  fixed
$\overline\tau>0,$ there exists $\rho=\rho(\overline\tau)>0$ such that
\be
\A(\tau) \subseteq  K^{\rho}
\qquad\ \  \forall~\tau \geq \overline\tau\,.\label{t2}
\ee
Moreover, there exist $T>0$ and $\rho'<\rho(T)$, such that
%
%
{\setlength\arraycolsep{2pt}
\begin{eqnarray}
K^{\rho'}  &\subseteq & \A(\tau)
\qquad\, \forall~\tau > T\,.
\label{t1}
\end{eqnarray}}
%
%
\end{theorem}
\begin{remark}
\label{remark21}
Observe that, given $\vfi\in K^\rho$,
any map $x \mapsto w_i(\vfi(x))$,
$i\in \{1,\dots,n\}$,
is essentially bounded and has finite total increasing
variation 
on subsets of $[a,\, b]$ bounded away
from the end points $a,\, b$.
Hence, any map $x \mapsto w_i(\vfi(x))$,
$i\in \{1,\dots,n\}$,
has also finite total variation
on such sets and, in particular,
it admits
left and right limits in any point $x\in ]a, \, b[$.
Moreover, since
an element $\vfi$ of $K^\rho$ is 
defined up to $\lu$ equivalence,
we may always assume that there is
a right continuous representative of 
$w_i(\vfi)$, $i\in \{1,\dots,n\}$,
that satisfies the inequalities 
appearing in the definition of $K^\rho$.
\end{remark}

\begin{theorem} \label{theorem2}
Under the same assumptions of Theorem 1, the set $\A(T)$ is 
a compact subset of~$\lu ([a,b],\,\Gamma)$ for
each $T>0.$ 
\end{theorem}

Indeed, we will prove in Section~\ref{section34}
that the compactness of the attainable set $\A(T)$ holds 
even in the case where $\A(T)$ is defined
as the set of all configurations that can be reached
at time $T$ only
by solutions of the mixed problem for (\ref{311})
that admit a strong $\lu$ trace at the boundaries $x=a, \, x=b$
(as the ones generated by a front tracking algorithm).

\vsp
\subsection{Existence and uniqueness of solutions}
\label{subsec:existence}
\indent


We describe here a front tracking algorithm 
that generates approximate solutions to (\ref{311})
on the strip $[a,b]\times \R^+$
continuously depending on the initial and boundary data,
which represents
a natural extension of \cite{A-G, brgoa1}.
Fix an integer $\nu\geq 1$ and consider the discrete set of points in $\Gamma$
whose coordinates are integer multiples of $2^{-\nu}$:
\be
\Gamma^\nu\doteq \Big\{ u\in \Gamma\,;~ w_i(u)\in 2^{-\nu}\Z,\quad 
i=1,\ldots,n\Big\}.
\label{gammappr}
\ee
Moreover, consider the domain
\be
\D^\nu\!\doteq\!
\Big\{(u, u', u'');\quad
 u\in \li([a,b],\Gamma^\nu),\ u', u''\in \li(\R^+,\Gamma^\nu),\ \
u,\, u', \, u''~~\hbox{are piecewise
constant}\Big\}.
\label{332} 
\ee
On $\D^\nu$ we now construct a flow map $E^\nu$ whose trajectories
are front tracking approximate solutions of (\ref{311}).
To this end, we first describe how to solve a Riemann problem
with left and right initial states $u^L,u^R\in \Gamma^\nu$.  
In Riemann coordinates, assume that
$$
w(u^L)\doteq w^L=(w^L_1,\ldots,w_n^L),\qquad 
w(u^R)\doteq w^R=(w^R_1,\ldots,w_n^R).
$$
Consider the intermediate states 
\be
z^0\doteq u^L,\qquad\ldots,\qquad \quad
z^i\doteq u(w^R_1,\ldots,w^R_i,w^L_{i+1},\ldots,w^L_n),\qquad\ldots,\qquad
z^n \doteq u^R.\label{333}
\ee
The solution to the Riemann problem $(u^L, u^R)$ is constructed
by piecing together the solutions to the simple
Riemann problems $(z^{i-1}, z^i)$, $i=1, \dots ,n$.
If $w_i^R<w_i^L$, the solution 
of the Riemann problems $(z^{i-1}, z^i)$
will contain a single $i$-shock,
connecting the states $z^{i-1}, \,z^i,$ and traveling
with the Rankine-Hugoniot speed $\lambda_i(z^{i-1},\,z^i)$.
Here and in the sequel, 
by $\lambda_i(u,u')$ we denote the $i$-th eigenvalue of the
averaged matrix
\be
A(u,u')\doteq \int_0^1 
D f\big(\theta u+(1-\theta)u'\big)\,d\theta.\label{334}
\ee
If $w_i^R>w_i^L$, the exact solution of the Riemann problem 
$(z^{i-1}, z^i)$ would contain
a centered rarefaction wave.  This is approximated by a rarefaction
fan as follows.   If $w_i^R=w_i^L+ p_i \,2^{-\nu}$ we insert the states
\be
z^{i,\ell}\doteq(w^R_1,\ldots,w^L_i+\ell \,2^{-\nu},w^L_{i+1},\ldots,w^L_n),
\qquad\qquad\ell=0,\ldots,p_i,\label{335}
\ee
so that $z^{i,0}=z^{i-1}$, $~z^{i,p_i}=z^i$.
Our front tracking solution will then contain
$p_i$ fronts of the $i$-th family, each connecting a couple of states
$z^{i,\ell-1}$, $z^{i,\ell}$ and traveling with speed
$\lambda_i\big(z^{i,\ell-1},~z^{i,\ell}\big)$.
\vsp

For any given triple  of (piecewise constant)
initial and boundary data $(\,\overline u,\,\widetilde u_a\,,\,\widetilde
u_b\,)\in\D^\nu$, the approximate solution 
$u(t,\cdot)\doteq E_t^\nu
(\,\overline u,\,\widetilde u_a\,,\,\widetilde u_b\,)$
is now constructed as follows. At time $t=0,$ 
for $a<x<b$ we solve the initial
Riemann problems
determined by the jumps in $\overline u$ according 
to the above procedure,
while 
at $x=a$ we construct the solution to the Riemann problem 
with left and right initial states  
$u^L=\widetilde u_a(0+)$, $u^R=\ov u(a+)$ and take its restriction to the 
interior of the domain $\Omega$. 
In the same way, at $x=b$, we take  the restriction to the 
interior of $\Omega$ of the solution to the Riemann problem 
with initial states $u^L=\ov u(b-)$, $u^R=\widetilde u_b(0+)$. 
This yields
a piecewise constant function with finitely many fronts, traveling with
constant speeds. The solution is then  prolonged up to the first
time $t_1$ at which one of the following events takes place:

\vspace{5pt}

a) two or more discontinuities interact in the interior of $\Omega$;

b) one or more discontinuities hit the boundary of $\Omega$;

c) the boundary data $\widetilde u_a$ has a jump;

d) the boundary data $\widetilde u_b$ has a jump.

\vsp

\n 
If the case a) occurs, we then solve
the resulting Riemann problems 
applying again the above procedure, while in the other three cases b)-c)-d)
we construct the solution to the Riemann problem 
with left and right initial states  
$u^L=\widetilde u_a(t_1+)$, $u^R=u(t_1,a+)$, or $u^L= u(t_1, b-)$,
$u^R=\widetilde u_b(t_1+)$,  and take its restriction to the  interior of the
domain $\Omega$. This determines the solution $u(t,\cdot)$ until the time
$t_2>t_1$ where one of the events a),b),c) again takes place,
etc$\ldots$  
Notice that at any time
where case b) occurs but c) or d) do not take place, no new wave is
generated. Therefore,  waves entering
the domain $\Omega$ at the boundaries $x=a, \, x=b$
are produced only by the
jumps of the boundary data $\widetilde u_a,\> \widetilde u_b$.
\vspace{5pt}

As in \cite{A-G, brgoa1}, one checks that 
the approximate solution $u$
constructed with this algorithm is 
well defined for all times $t\geq 0$.  
Indeed, the following properties hold.
\begin{itemize}
\item[-] The total variation of $u(t,\cdot),$
measured w.r.t. the Riemann coordinates $w_1(t,\cdot),\dots,$
$w_n(t,\cdot)$,
is non-increasing in time.
\item[-] The number of wave-fronts in $u(t,\cdot)$ is non-increasing
at each interaction. Hence, the total number of wave-fronts in 
$u(t,\cdot)$  remains finite.
\end{itemize}
It is then possible to define a flow map
\be
\p \mapsto 
E_t^\nu \p,
\qquad\quad 
\p\doteq(\,\overline u,\,\widetilde u_a\,,\,\widetilde u_b\,)\in \D^{\nu},
\quad t\ge 0
\label{flowappr}
\ee
of approximate solutions of (\ref{311}). 
By construction, each trajectory
$t \mapsto 
E_t^\nu \p$
is a weak  solution of (\ref{311}) (because all
fronts of $u(t,\cdot)\doteq E_t^\nu \p$
satisfy the Rankine-Hugoniot conditions),
but may contain discontinuities that do not
satisfy the usual Lax stability conditions
(due to the presence of rarefaction fronts).
On the other hand, one can verify as in~\cite[Lemma~4.4]{A-G}
that, due to genuine nonlinearity, the amount of positive
waves in $u(t,\cdot),$
measured w.r.t. the Riemann coordinates 
$w_1(t,\cdot),\dots,w_n(t,\cdot)$,
decays in time and in space.
Hence, for a.e. $a<x<y<b,$ one obtains the
Oleinik type estimates 
\be
{\setlength\arraycolsep{0pt}
\begin{array}{ll}
w_i(t,y)-w_i(t,x) ~&\leq ~C\cdot 
\displaystyle{\bigg\{{{y-x}\over
{t}}+\log\bigg(\frac{y-b}{x-b}\bigg)\bigg\}
+N_\nu\, 2^{-\nu}}
\quad\ \textrm{if}\quad \ \  i\in\{1,\ldots,p\}, 
\\
\noalign{\medskip}
w_i(t,y)-w_i(t,x) ~&\leq ~
C\cdot
\displaystyle{\bigg\{{{y-x}\over 
{t}}+\log\bigg(\frac{y-a}{x-a}\bigg)\bigg\}
+N_\nu\, 2^{-\nu}}
\quad\ \textrm{if}\quad \ \ i\in\{p+1,\ldots,n\},
\end{array}
}
\label{oleinappr1}
\ee
%
%
where
$N_\nu$ denotes the maximum number of shocks
of each family present in the initial data $\overline u$,
and in the boundary data $\widetilde u_a\,,\,\widetilde u_b$.
Similarly, one can check that along the $x$-sections,
for a.e. $0<\tau_1<\tau_2$,
there holds
\be
{\setlength\arraycolsep{0pt}
\begin{array}{ll}
w_i(\tau_2,x)-w_i(\tau_1,x) ~&\leq ~C\cdot
\displaystyle{\bigg\{{{\tau_2-\tau_1}\over {x-b}}+
\log\left(\frac{\tau_2}{\tau_1}\right)\!\!\bigg\}
+N_\nu\, 2^{-\nu}}
\quad\ \textrm{if}\quad\ \ i\in\{1,\ldots,p\}, 
\\
\noalign{\medskip}
w_i(\tau_2,x)-w_i(\tau_1,x) ~&\leq
~C\cdot
\displaystyle{\bigg\{{{\tau_2-\tau_1}\over {x-a}}+
\log\left(\frac{\tau_2}{\tau_1}\right)\!\!\bigg\}
+N_\nu\, 2^{-\nu}} 
\quad\ \textrm{if}\quad \ \ i\in\{p+1,\ldots,n\}.
\end{array}}
\label{oleinappr2}
\ee
\begin{remark}
\label{remark24}
Observe that, if $u(t,x)$ is a front tracking solution
of the Cauchy problem for~(\ref{311})
(with initial data $\overline u(x)\doteq u(0,x)$)
constructed by the algorithm  in \cite{brgoa1}
on the upper half \linebreak plane~$\R^+\times \R$,
then the restriction of $u(t, \cdot)$ to the interval
$[a,b]$ coincides with the front tracking solution
$E_t^\nu
(\,\overline u,\,\widetilde u_a\,,\,\widetilde u_b\,)$ 
of the mixed problem for (\ref{311}),
with
boundary data 
$
\widetilde u_a(t)\doteq u(t,a)$, \,
$\widetilde u_b(t) \doteq u(t,b)\,.
$
\end{remark}
\vspace{5pt}
\vspace{5pt}

As $\nu \to \infty$, the domains $\D^{\nu}$
become dense in 
\be \D \doteq
\Big\{ (\,\overline u,\,\widetilde u_a\,\,\widetilde u_b\,)~;~
\overline u\in\li ([a,b],\Gamma),\>\> \widetilde u_a, \,\widetilde u_b\in\li
(\R^+,\Gamma) \Big\}. 
\label{dominio}
\ee
Thus, following the same technique adopted in~\cite{A-G},
one can define a flow
map $E_t$ on $\D$ as a suitable limit of the flows $E_t^\nu$
in (\ref{flowappr}),
that depends Lispschitz continuously on the initial and boundary
data.
Namely, the following holds.
\begin{theorem} \label{theorem31}
Let  (\ref{311}) be a system of Temple class
with all characteristic fields genuinely nonlinear, and assume
that the strict hyperbolicity condition ${\bf (SH1)}$ holds.
Then, there exists a  continuous map 
\be
(t,\, \overline u,\,\widetilde u_a\,,\,\widetilde u_b\,)
\mapsto E_t  (\,\overline u,\,\widetilde u_a\,,\,\widetilde u_b\,)
\qquad\quad t\geq 0, \ \
(\,\overline u,\,\widetilde u_a\,,\,\widetilde u_b\,)\in \D\,,
\label{flowmap}
\ee
and some
constant $C>0$ depending only on the system (\ref{311})  and on the domain
$\Gamma,$ so that, \linebreak
\pagebreak

\n
for every fixed  $0<\delta<(b-a)/2$, and for all 
$\p_1\doteq
(\,\overline u,\,\widetilde u_a\,,\,\widetilde u_b\,), \,
\p_2\doteq(\,\overline v, \,\widetilde v_a\,, \,\widetilde v_b\,)\in\D,$
letting $L_t\doteq L_t(\delta)=C(1+\log(t/\delta))$, 
there holds
\be
{\setlength\arraycolsep{2pt}
\begin{array}{ll}
&\big\|E_t \p_1-
E_t \p_2
\big\|_{{\LL}^{^{\!1}}([a+\delta,\,b-\delta])} \leq 
L_t\cdot 
\Big\{\big\|\,\ov u-\ov v\,\big\|_{{\LL}^{^{\!1}} ([a,b])} +
\big\|\widetilde u_a-\widetilde v_a\big\|_{{\LL}^{^{\!1}}([0,\,t])}+
%
\big\|\widetilde u_b-\widetilde v_b\big\|_{{\bf L}^{^{\!1}}([0,\,t])} 
\Big\}
%
\\
\noalign{\smallskip}
\end{array}
}
\label{3217}
\ee
for all $t\geq \delta$.
Moreover, the map $(t,x)\mapsto E_t (\,\overline
u,\,\widetilde u_a\,,\,\widetilde u_b\,)(x)$  yields an entropy weak solution
 (in the sense of Definition \ref{definition32}) to the  initial-boundary
value problem (\ref{311})-(\ref{314}) on $\Omega,$ that admits strong
$\ilu$~traces at the boundaries $x=a$ and $x=b$, i.e. there exist two
measurable maps $\psi_a,\, \psi_b: \R^+ \mapsto U$  such that 
\be
{\setlength\arraycolsep{2pt}
\begin{array}{ll}
\displaystyle{
\lim_{x\to a^+} \int_0^\tau \big|E_t (\,\overline u,\,\widetilde
u_a\,,\,\widetilde u_b\,)(x)-  \psi_a(t)\big|~dt} 
&=0,
\\
\noalign{\medskip}
\displaystyle{
\lim_{x\to b^-} \int_0^\tau
\big|E_t (\,\overline u,\,\widetilde u_a\,,\,\widetilde u_b\,)(x)- 
\psi_b(t)\big|~dt }
&=0, 
\end{array}
}
\qquad\ \forall~\tau\geq 0.
\label{strongtrace1}
\ee
\end{theorem}
\vsp
The proof of Theorem~\ref{theorem31}
can be obtained with entirely similar arguments
to those used to establish 
\cite[Theorem~2.1]{A-G}, where 
a continuous flow of solutions to (\ref{311}) is
constructed 
in the case of a mixed problem 
on the
quarter of plane $\{(t,x)\in \R^2~;~t\geq 0,\, x\geq 0\}$,
with a single boundary at $x=0$.
\vsp

Concerning uniqueness,
with the same arguments in~\cite{A-G} one obtains
the following
result which is the extension of \cite[Theorem~2.2]{A-G}
to the present case of a domain $\Omega$ with
two boundaries at $x=a$ and at $x=b$.
\vspace{5pt}

\begin{theorem} \label{theorem32}
Let  (\ref{311}) be a system of Temple class 
satisfying the same assumptions as in 
Theorem \ref{theorem31}. Let $u=u(t,x)$ be an entropy weak
solution to the mixed problem (\ref{311})-(\ref{314}) 
on the region $\Omega_T\doteq[0,\,T]\times [a,b]$
(in the sense of Definition \ref{definition32}).
Assume that the following conditions hold.
\begin{itemize}
\item[(i)] The map $(t,x) \to \big(u(t, \cdot),\,u(\cdot, x)\big)$
takes values within the domain 
\be
\D_T\doteq\Big\{(\,\overline u,\,\widetilde u_a,\,\widetilde
u_b\,)~;~\overline u\in\li ([a,b]
,\Gamma),\,  \widetilde u_a,\,\widetilde u_b \in\li ([0,T],\Gamma)
\Big\}.
\label{domainT}
\ee
\item[(ii)] There holds
\be
\label{qcinitial}
\displaystyle{\ess \sup_{t\to 0^+}} \int_a^b
\big|u(t,x)-\overline u(x)\big| \,dx=0\,.
\ee
\item[(iii)] 
There holds
{\setlength\arraycolsep{2pt}
\begin{eqnarray}
\label{qcboundary1}
\displaystyle{\ess \sup_{x\to a^+}}\int_0^T
\big|w_i\big(u(t,x)\big)-w_i(\widetilde u_a(t))\big| \, dt&=&0
\qquad\quad\forall~i=p+1,\dots,n, \qquad\qquad
\\
\noalign{\medskip}
\label{qcboundary2}
\displaystyle{\ess \sup_{x\to b^-}}\int_0^T
\big|w_i\big(u(t,x)\big)-w_i(\widetilde u_b(t))\big| \, dt&=&0
\qquad\quad\forall~i=1,\dots,p. \qquad\qquad
\end{eqnarray}}
\end{itemize}

\noindent
Then, $u$ coincides with the corresponding trajectory of the
flow map $E_t$ provided by Theorem~\ref{theorem31},
namely one has
\be
u(t, \cdot) = E_t(\,\ov u,\,\widetilde u_a\,,\,\widetilde u_b)( \cdot),\qquad\ 
\forall~0\leq t \leq T.
\label{3218}
\ee
\end{theorem}
\vsp

The next result shows that
the conditions (\ref{qcinitial})-(\ref{qcboundary2})
are certainly satisfied
by entropy weak  solutions
to the mixed problem
(\ref{311})-(\ref{314}) obtained as limit  
of front tracking approximations.

\begin{theorem}
\label{theorem25}
Let  (\ref{311}) be a system of Temple class
satisfying the same assumptions as \linebreak in
Theorem \ref{theorem31}.
Consider a sequence
$u^\nu(t, \cdot) : [a, b] \mapsto \Gamma^\nu$
of wave-front tracking
approximate solutions of the mixed problem for~(\ref{311})
(constructed with the above algorithm)
that converges in~$\lu$, as $\nu \to \infty$, to some function
$u(t, \cdot) : [a, b] \mapsto \Gamma$,
for every $t\in [0,T]$.
Then, there exist the right limit
at $x=a$, and the left limit at $x=b$,
of the  map $x \to u(t,x)$ for every $t\in [0, T]$,
and the right limit at $t=0$ of the map $t \to u(t,x)$
for every $x\in [a, b]$.
Moreover, there is a countable set~$\N\subset \R$
such that $u(t,a)=u(t, a+), \, u(t, b)
=u(t, b-)
$
for all  $t\in [0,T]\setminus \N$, and 
$u(0,x)=u(0+, x)
$ for all $x\in [a,b]\setminus \N$,
and, setting 
$\overline u \doteq u(0, \cdot)$, $\widetilde u_a \doteq
u(\cdot,a), \, \widetilde u_b\doteq u(\cdot,b)$,
there holds (\ref{3218}).
\end{theorem}
\begin{remark}
\label{remark25}
It was shown in~\cite[Lemma~2.1]{A-G} that an  alternative way to prove
the essential  \linebreak 
limits~(\ref{qcboundary1})-(\ref{qcboundary2}),
is to employ the distributional
entropy inequalities associated to the
``boundary entropy pairs'' for~(\ref{311}),
introduced by
G.-Q.~Chen and H.~Frid in \cite{chfr1, chfr2}.
However, in order to apply~\cite[Lemma~2.1]{A-G}
to a function $u$ obtained as a limit of approximate
solutions~$u^\nu$, it is necessary to know
the $\lu$ convergence of the sequence of the 
corresponding boundary 
data~$\widetilde u^\nu_a, \, \widetilde u^\nu_b$. 
Instead, the result provided 
here by Theorem~\ref{theorem25} allows to derive 
the limits~(\ref{qcboundary1})-(\ref{qcboundary2})
requiring only the $\lu$ convergence of the sequence
of the approximate solutions
$u^\nu(t, \cdot)$, for all $t$. This property
will be crucial to establish the main result
of the paper stated in Theorems~\ref{theorem1}-\ref{theorem2}.
\end{remark}
\vsp

In order to prove Theorem~\ref{theorem25}, we will
show in the next section that, for Temple systems,
solutions of the mixed problem
(\ref{311})-(\ref{314}) with possibly unbounded 
variation 
enjoy the same
regularity property
(of being continuous outside a countable
number of Lipschitz curves) possessed
by solutions
with small total variation
of  a general 
system, thus extending the regularity
results obtained under the smallness assumption
of the total variation
by DiPerna~\cite{dip}
and Liu~\cite{liu}
(for solutions constructed by the Glimm scheme) and
by Bressan and LeFloch
\cite{brfiocco} (for solutions
generated by a front tracking algorithm).

\begin{prop}
\label{theorem26}
In the same setting as Theorem~\ref{theorem25},
consider a sequence
$u^\nu(t, \cdot) : [a, b] \mapsto \Gamma^\nu$
of wave-front tracking
approximate solutions of the mixed problem for~(\ref{311})
(constructed with the above algorithm)
that converges in $\lu$, as $\nu \to \infty$, to some function
$u(t, \cdot) : [a, b] \mapsto \Gamma$,
for every $t\in [0,T]$.
%
%
%
Then,
%
%
there exist a countable set of interaction points
$\Theta\doteq \big\{ (\tau_l,\, x_l ); ~l\in
\NN\big\}\subset \Omega_T\doteq [0,\,T]\times [a,\, b]$,
and a countable family of Lipschitz continuous shock curves
$\Upsilon \doteq \big\{ x=y_m(t);~t\in\, ]r_m,\,s_m[,\>m\in
\NN\big\}$, such that the following hold.
\begin{itemize}
\item[(i)]
For each $m\in\NN$, and
for any $\tau\in \,]r_m,\,s_m[$ with $(\tau, y_m(\tau))\not\in
\Theta$, there exist the derivative
$\dot y_m(\tau)$ and the left and right  limits
\be \lim\limits_{(s,y)\rightarrow (\tau, y_m(\tau)),\, y< y_m(\tau)} u(s,y)
\doteq u^-,\qquad
\lim\limits_{(s,y)\rightarrow (\tau, y_m(\tau)),\, y> y_m(\tau)} u(s,y)
\doteq u^+\,.
\label{limiti}\ee
Moreover, these limits satisfy the Rankine Hugoniot relations
\be \dot y_m(\tau) \cdot (u^+-u^-)= f(u^+)- f(u^-)
\label{RH}\ee
and, for some $i\in \{1,...,n\}$, there hold
the Lax entropy inequalities
\be \lambda_i(u^+)<\dot y_m(t) < \lambda_i(u^-)\,. \label{LAX}\ee
\item[(ii)]
The map $u$ is continuous outside the set $\Theta \cup \Upsilon.$
\end{itemize}

\end{prop}
\vsp

%
%
\section{Preliminary results} \label{section33}
\setcounter{equation}{0}

In this section we first
provide some
estimates on the distance between two rarefaction fronts
of a front tracking solution (constructed by the algorithm described
in Section~\ref{subsec:existence})
similar \linebreak to~\cite[Lemma~4]{brgoa1},  \cite[Prop.~4.5]{bianc}.
We next show how to
approximate the profile $u(t,\cdot)$ of a solution of the mixed problem
(\ref{311})-(\ref{314}), with a function taking values in the
discrete set~$\Gamma^\nu$ defined at~(\ref{gammappr}),
which enjoys the same type of estimates on the positive
waves as~$u(t,\cdot)$.
We conclude the section
establishing the regularity result stated in 
Proposition~\ref{theorem26} on the
global structure of solutions to the mixed problem
for (\ref{311}), which in turn yields Theorem~\ref{theorem25}.

%
\begin{lemma} \label{lemma1}
There exists some constant $C_1 \,>0$
depending only on the system (\ref{311})
such that the following holds.
Consider a front tracking solution $u(t,x)$
with values in   $\Gamma^\nu$,
constructed by the algorithm  of Section~\ref{subsec:existence}
on the region $[\tau,\tau']\times [a,b]$.
Then, given
any two adjacent rarefaction fronts of $u$ located at
$x(t)\leq y(t)$, $t\in [\tau, \, \tau']$, and
belonging to the same family,
there holds
%
%
\be
\big| y(\tau') - x(\tau') \big| 
\le \,\big| y(\tau) - x(\tau) \big| + C_1 (\tau'- \tau)\,
2^{-\nu}\,.
\label{distraref}
\ee
\end{lemma}
\vspace{5pt}
\n\textsc{Proof.}
Consider two 
adjacent rarefaction fronts of the $k$-th 
family
$x(t) \leq y(t)$, $t\in[\tau,\,\tau']$,
 and let  $\tau_1 <...< \tau_N$ 
be the interaction times of 
$x(t)$ in the interval $[\tau, \tau']$. Set
$\tau_0 \doteq \tau$, $\tau_{N+1} \doteq \tau'$, and
fix  $\alpha\in\{0,...,N\}$.
Let $t \to z(t; s,x)$ be 
the characteristic curve of the $k$-th
family starting at $(s, x),$ i.e.
the solution to the ODE
$$
\dot z = \lambda_k (u(t,z)), \qquad z(s; s, x)= x.
$$
Notice that, although the above ODE has discontinuous
right hand-side (because of the discontinuities
in the front tacking solution
$u$),  its solution $z(\cdot; s,x)$ is unique
and depends Lipschitz continuously on the initial data $x$
since it crosses only a finite number of 
jumps~(see~\cite{br5}).
Choose $t_0<t_1<{\tau}_{\alpha+1}$ so that 
the characteristic curve
$z(\cdot;t_0, x(t_0))$ 
does not cross any
wave-front of the other families in the interval  
$[t_0,\, t_1]$, and then, by induction,
define a sequence
of times  
$\{t_i\}_{i\in\Z}\,\subset ~]\tau_\alpha, \tau_{\alpha +1}[\,$ 
so that
\be
\begin{array}{ll}
&\tau_\alpha < t_{-i-1} < t_{-i} 
\le t_0 \le t_{i}< t_{i+1}<\tau_{\alpha +1}, 
\quad i\in\NN,
\\
\noalign{\medskip}
&\hskip 0.5in \lim\limits_{i\rightarrow -\infty} t_i = \tau_{\alpha}, \qquad 
\lim\limits_{i\rightarrow +\infty} t_i = \tau_{\alpha+1},
\end{array}
\label{sequence}
\ee
with the properties that
the characteristic curve of the $k$-th
family starting at $(t_i, x(t_i)),$
does not cross any
wave-front of the other families in the interval  
$[t_i,\, t_{i+1}]$, for each $i\in\Z.$
Thus, setting
$$
u^+_i \doteq u(t_i,\, x(t_i)+),
\qquad u^-_i \doteq u(t_i,\, x(t_i)-),
$$
and observing that, by construction, one has 
$|w(u^+_i) - w(u^-_i)| < 2^{-\nu}$,
we derive
\be
{\setlength\arraycolsep{2pt}
\begin{array}{ll}
\big|z(t_{i+1};\,t_i,\, x(t_i))-x(t_{i+1})\big|
& \leq (t_{i+1}-t_i)\cdot
\big|\lambda_k(u^+_i)-\lambda_k(u^+_i,\,u^-_i)\big|
\\
\noalign{\medskip}
& \leq c\cdot(t_{i+1}-t_i)\cdot\big|w(u^+_i) - w(u^-_i)\big|
\\
\noalign{\medskip}
& \leq c \cdot(t_{i+1}-t_i)\cdot 
2^{-\nu}
\end{array}}
\label{stefano1}
\ee
\vspace{-8pt}

\n
for some constant $c>0$ depending only on the system.
Relying on (\ref{stefano1}), and
since $z(\tau';\, t_{i+1},\,x)$  depends Lipschits continuously 
on the initial data $x$,
we deduce
that there exists some other constant $c'>0$, depending only
on the system and on the set $\Gamma$, so that
there holds
\be
{\setlength\arraycolsep{2pt}
\begin{array}{ll}
\big|z(\tau';\,t_i,\, x(t_i))-z(\tau';\,t_{i+1},\, x(t_{i+1}))\big|
& \leq c'\cdot\big|z(t_{i+1};\,t_i,\, x(t_i))-x(t_{i+1})\big|
%
\\
\noalign{\medskip}
& \leq c'\cdot c \cdot(t_{i+1}-t_i)\cdot 
2^{-\nu}
\end{array}}
\label{stefano2}
\ee
for any $i\in\Z.$ 
Thus, by (\ref{sequence}),
and thanks to (\ref{stefano2}), we obtain
\be
{\setlength\arraycolsep{2pt}
\begin{array}{ll}
\big|z(\tau';\,\tau_\alpha,\, x(\tau_\alpha))-z(\tau';\,\tau_{\alpha+1},\, 
x(\tau_\alpha))\big|
& \leq 
\displaystyle{
\sum\limits_{i\in\Z}\big|z(\tau';\,t_i,x(t_i)) -
z(\tau';\,t_{i+1},\,,x(t_{i+1}))\big|}
%
\\
\noalign{\medskip}
& \leq c'\cdot c \cdot(\tau_{\alpha+1}-\tau_{\alpha})\cdot 
2^{-\nu}\,.
\end{array}}
\label{stefano3}
\ee
Repeating this computation for every interval
$]\tau_\alpha, \tau_{\alpha +1}[$,
$\alpha \in \{0,...,N\},$ we get
\be
{\setlength\arraycolsep{2pt}
\begin{array}{ll}
\big|z(\tau';\,\tau,\, x(\tau))-x(\tau')\big|
& \leq 
\displaystyle{
\sum\limits_{\alpha=0}^N
\big|z(\tau';\,\tau_\alpha,\, x(\tau_\alpha))-z(\tau';\,\tau_{\alpha+1},\, 
x(\tau_\alpha))\big|
}
%
\\
\noalign{\medskip}
& \leq c'\cdot c \cdot(\tau'-\tau)\cdot 
2^{-\nu}\,.
\end{array}}
\label{stefano4}
\ee
Clearly, one obtains the same type of estimate as (\ref{stefano4})
for the other rarefaction front $y(t)$, i.e. there holds
\be
\big|z(\tau';\,\tau,\, y(\tau))-y(\tau')\big|
\leq c'\cdot c \cdot(\tau'-\tau)\cdot 
2^{-\nu}\,.
\label{stefano5}
\ee
On the other hand, by (\ref{uppbound}),
we have
\be
\big|z(\tau';\,\tau,\, x(\tau))-
z(\tau';\,\tau,\, y(\tau))\big|\leq
\big|x(\tau)-y(\tau)\big|+ 2\,\lambda^{\rm max}\cdot (\tau'-\tau).
\label{stefano6}
\ee
Thus, (\ref{stefano4})-(\ref{stefano6}) together
yield (\ref{distraref}),
concluding the proof.
\hfill$\Box$
\vs

In the following, in connection with any 
(right continuous) piecewise constant 
map \linebreak $\psi:[a,b]\mapsto 2^{-\nu}\,\Z$,
we will let $\pi(\psi)=\{x_0=a<x_1<\cdots <x_{\overline \ell}=b\}$ 
denote the partition of~$[a,b]$ induced by
$\psi$, in the sense that $\psi(x)$ is constant 
on every interval $[x_{\ell}, \, x_{\ell+1}[$\,,
$0 \leq \ell <\overline \ell$.
Then, given $\rho>0$, for any $\nu\ge1,$ 
consider the set of piecewise constant maps
\be 
\begin{array}{ll}
\!\!\!\!\!\!\!\!\!\!
&K^\rho_{\nu}\doteq\!
\left\{ \vfi:[a,b]\mapsto \Gamma^{\nu};
\!\!\begin{array}{l}
\displaystyle{\frac{w_i(\vfi(x_k))-w_i(\vfi(x_h))}{x_k-x_h}}
\leq \displaystyle{\frac{5\rho}{x_h-a}}\qquad 
\left\{
\!\!\!\!\!\!\!\!
\begin{array}{ll}
&\!\hbox{for}
\begin{array}{ll}
&a<x_h<x_k<b, 
\\
\noalign{\smallskip}
&x_h,\,x_k\in \pi(w_i\circ\vfi),
\end{array}
\\
\noalign{\medskip}
&\!\hbox{if} \qquad i\in \{p+1,...,n\}
\end{array}
\right.
\\
\noalign{\bigskip}
\displaystyle\frac{w_i(\vfi(x_k))-w_i(\vfi(x_h))}{x_k-x_h}\leq 
\displaystyle{\frac{5\rho}{b-x_k}}\qquad 
\left\{
\!\!\!\!\!\!\!\!
\begin{array}{ll}
&\!\hbox{for}
\begin{array}{ll}
&a<x_h<x_k<b, 
\\
\noalign{\smallskip}
&x_h,\,x_k\in \pi(w_i\circ\vfi),
\end{array}
\\
\noalign{\medskip}
&\!\hbox{if} \qquad i\in \{1,...,p\}
\end{array}
\right.
\end{array}
\!\!\!\!\!\!\right\}.
\\
\noalign{\smallskip}
\end{array} 
\label{Kro}  
\ee
The next lemma shows that we can approximate
in $\lu$ any map $\vfi\in K^\rho$ with 
a piecewise constant function~$\vfi_\nu\in K^\rho_\nu$.
\vskip 5pt

\begin{lemma} \label{lemma2}
For any given $\vfi\in K^\rho$,
there exists a sequence 
of right continuous 
maps~$\vfi_{\nu}\in K^{\rho}_{\nu}$, $\nu \geq 1$,
such that:
\begin{itemize}
\itemindent 8pt
\item[{\it a)}] for every $i\in \{1,\dots,n\}$,
and for any $x_h\in \pi(w_i\circ\vfi_\nu)$, there holds
\be
w_i(\vfi_\nu(x_{h+1}))>w_i(\vfi_\nu(x_h))
\qquad \Longrightarrow
\qquad 
w_i(\vfi(x_{h+1}))=w_i(\vfi(x_h))+2^{-\nu}\,;
\label{rarefcond}
\ee
\item[{\it b)}] there holds
\be
\vfi_{\nu}\rightarrow \vfi\qquad {\rm in }\qquad \lu([a,b]). 
\label{approxconv} 
\ee
\end{itemize}
\end{lemma}
\vspace{5pt}

\noindent
{\bf 1.}
First observe that,
by Remark~\ref{remark21}, any map $x \mapsto w_i(\vfi(x))$, 
$i \in \{1,\dots,n\}$
has finite total variation on $[a+\varepsilon, \, b-\varepsilon],\>\varepsilon>0$.
Hence, we may assume that 
$w_i(\vfi(\cdot))$ admits left and right  limits
in any \linebreak point $x\in ]a, \, b[$,
and that $w_i(\vfi(x))=w_i(\vfi(x^+))\doteq 
\lim_{\xi\to x^+} w_i(\vfi(\xi))$,  for
all~$i \in \{1,\dots,n\}$.
Let $\{y_{i,m}~;~m\in\NN\}$ be the countable set
of discontinuities  of $w_i(\vfi(\cdot))$, 
$i \in \{1,\dots,n\}$.
Then, 
we can find a partition 
$\xi_{i,m}^1= y_{i,m}<\xi_{i,m}^2<\cdots <\xi_{i,m}^{\ell_{i,m}}
=y_{i,m'}$  of each interval $[y_{i,m}, \, y_{i,m'}[$
where $x \mapsto w_i(\vfi(x))$ is continuous, so that:
\begin{itemize}
\itemindent 8pt
\item[\textsl{i)}]
\ for every $1 <\ell<\ell_{i,m}$ there holds
\be
w_i(\vfi(\xi_{i,m}^\ell))\in 2^{-\nu}\,\Z \,;
\label{appr1}
\ee
\item[\textsl{ii)}]
\ for every $1 \leq\ell<\ell_{i,m}$ one has
\be
\big|w_i(\vfi(x))-
w_i(\vfi(\xi_{i,m}^\ell))
\big|
\leq 2^{-\nu}
\qquad\quad\forall~x\in[\xi_{i,m}^\ell,\,\xi_{i,m}^{\ell+1}[\,.
\label{appr2}
\ee
%
%
\end{itemize}
Notice that  the Oleinik type conditions
stated in the definition of $K^\rho$ imply that,
at any discontinuity point $y_{i,m}$ of $w_i(\vfi(\cdot))$, one has
\be
\lim_{\xi\to y_{i,m}^-} w_i(\vfi(\xi)) >
w_i(\vfi(y_{i,m}))\,.
\label{downdisc}
\ee
%
\vsp

\noindent
{\bf 2.} Let $\vfi_\nu:[a,b]\mapsto \Gamma^{\nu}$
be the piecewise constant, right continuous
map defined by setting, for every 
$i\in\{1,\dots,n\}$, and for any interval
$[y_{i,m}, \, y_{i,m'}[$ where 
$w_i(\vfi(\cdot))$ is continuous,
\be
 w_i(\vfi_\nu(x))\!\doteq\!
\left\{
{\setlength\arraycolsep{2pt}
\!\!\begin{array}{ll}
2^{-\nu}\lfloor 
2^{\nu}\,w_i(\vfi(\xi_{i,m}^1)) \rfloor \quad\ &
\hbox{if}\quad\ 
\left\{
\!\!\!\begin{array}{ll}
&x\in [\xi_{i,m}^1, \, \xi_{i,m}^2[\,,
\quad \hbox{and}
\\
\noalign{\smallskip}
&w_i(\vfi_\nu(\xi_{i,m}^1))\leq 
2^{-\nu}
\big(\lfloor 
2^{\nu}\,w_i(\vfi(\xi_{i,m}^1)) \rfloor\!+
2^{-1}\big),
\end{array}
\right.
\\
\noalign{\medskip}
2^{-\nu}
\big(
\lfloor 
2^{\nu}\,w_i(\vfi(\xi_{i,m}^1)) \rfloor\!+\!1\big) \quad\ &
\hbox{if}\quad\ 
\left\{
\!\!\!\begin{array}{ll}
&x\in [\xi_{i,m}^1, \, \xi_{i,m}^2[\,,
\quad \hbox{and}
\\
\noalign{\smallskip}
&w_i(\vfi_\nu(\xi_{i,m}^1))>
2^{-\nu}
\big(\lfloor 
2^{\nu}\,w_i(\vfi(\xi_{i,m}^1)) \rfloor\!+
2^{-1}\big),
\end{array}
\right.
\\
\noalign{\medskip}
w_i(\vfi(\xi_{i,m}^\ell)) \quad\ &
\hbox{if}\quad\ \ \  x\in [\xi_{i,m}^\ell, \, \xi_{i,m}^{\ell+1}[\,,
\qquad 1<\ell<\ell_{i,m}\,,
\end{array}
}
\right.
\ee
where $\lfloor \cdot \rfloor$ denotes the integer part.
Notice that, by construction, and because of (\ref{appr1})-(\ref{appr2}),
(\ref{downdisc}),
the map 
$\vfi_\nu :[a,b]\mapsto \Gamma^{\nu}$ 
enjoys the following property
\be
\left.
\begin{array}{ll}
& w_i(\vfi_\nu(x_k))>w_i(\vfi_\nu(x_h))
\\
\noalign{\medskip}
&  x_h<x_k\in \pi(w_i\circ\vfi_\nu)
\end{array}
\right\}
\qquad \Longrightarrow
\qquad 
w_i(\vfi(x_k))>w_i(\vfi(x_h))+2^{-(\nu+1)}\,.
\label{appr4}
\ee
\vskip 2pt

\noindent
Therefore, since
$\vfi\in K^\rho$, relying on 
(\ref{appr2}),
(\ref{appr4}), 
we deduce that, for every
$w_i(\vfi_\nu(\cdot))$, \linebreak 
$i\in 
\{p+1,\dots,n\}$, and for any $x_h<x_k\in \pi(w_i\circ\vfi_\nu)$
such that $w_i(\vfi_\nu(x_k))>w_i(\vfi_\nu(x_h))$,
there holds
\be
{\setlength\arraycolsep{2pt}
\begin{array}{ll}
\displaystyle{\frac{w_i(\vfi_\nu(x_k))-w_i(\vfi_\nu(x_h))}
{x_k-x_h}}
&\leq
\displaystyle{\frac{w_i(\vfi(x_k))-w_i(\vfi(x_h))+
2^{-(\nu-1)}}{x_k-x_h}}
\\
\noalign{\medskip}
&\leq
\displaystyle{\frac{5\,\big(w_i(\vfi(x_k))-w_i(\vfi(x_h))\big)}
{x_k-x_h}}
\\
\noalign{\medskip}
&\leq
\displaystyle{\frac{5\,\rho}{x_h-a}}\,.
\end{array}
}
\label{oleynappr1}
\ee
Clearly, with the same computations, we can show that, for every
$w_i(\vfi_\nu(\cdot))$,
$i\in \{1,\dots,p\}$, and for any $x_h<x_k\in \pi(w_i\circ\vfi_\nu)$,
there holds
\be
\displaystyle{\frac{w_i(\vfi_\nu(x_k))-w_i(\vfi_\nu(x_h))}
{x_k-x_h}}
\leq
\displaystyle{\frac{5\,\rho}{b-x_k}}\,.
\label{oleynappr2}
\ee
The estimates (\ref{oleynappr1})-(\ref{oleynappr2}),
together, imply that $\vfi_\nu\in K^\rho_\nu$,
while (\ref{appr2})
yields (\ref{approxconv}).
On the other hand observe that,
by construction, and because of (\ref{downdisc}),
the map $\vfi_\nu$ satisfies condition (\ref{rarefcond}),
which completes the proof of the lemma.
\hfill$\Box$
\vsp

We  now
provide a further
estimate on the distance between two rarefaction fronts
of a front tracking solution that, at a fixed time $\tau$, attains a profile
belonging to the set (\ref{Kro}).
\begin{lemma} \label{lemma3}
Consider a front tracking solution $u(t,x)$
with values in   $\Gamma^\nu$,
 $\nu \geq 1$,constructed by the algorithm  of Section~\ref{subsec:existence}
on the region $[\tau,\tau']\times [a,b]$.
Assume that $u(\tau', \cdot)$ is right-continuous,
verifies condition a) of Lemma~\ref{lemma2}, and satisfies
\be
u(\tau', \cdot)\in K_\nu^{\rho'},
\qquad\quad
\rho'\doteq
\displaystyle{\frac{\lambda^{\min}}{6 \,C_1}}\,,
\label{oleyncond}
\ee
where
$\lambda^{\min}$, $C_1$, are
the minimum speed in (\ref{uppbound}),
and the constant
of Lemma~\ref{lemma1}.
Then, given
any two adjacent rarefaction fronts of $u$ located at
$x(t)\leq y(t)$, $t\in [\tau, \, \tau']$, and
belonging to the same family,
there holds
%
%
\be
 x(\tau)<y(\tau).
\label{distraref2}
\ee
\end{lemma}
\vspace{5pt}
\n\textsc{Proof.}
To fix the ideas, assume that $x(t)\leq y(t)$
are the locations of
two adjacent rarefaction fronts of the $k\in\{p+1, \dots ,n\}$ - th
family, and hence, by (\ref{321}), have positive speeds.
Observe that, by
condition a) of Lemma~\ref{lemma2}, one has
\be
w_k(u(\tau',\, y(\tau')))- w_k(u(\tau',\, x(\tau')))  =
2^{-\nu} \,.
\label{rarefcon2}
\ee
Moreover, since $u$ is a front tracking solution
constructed by the algorithm of Section~\ref{subsec:existence}
on the region $[\tau,\tau']\times [a,b]$, we can apply
Lemma~\ref{lemma1}. Thus, using (\ref{uppbound}), (\ref{distraref}),
(\ref{rarefcon2}),
and recalling the definition (\ref{Kro}) of $K_\nu^{\rho'}$,
we deduce
$$
{\setlength\arraycolsep{2pt}
\begin{array}{ll}
 y(\tau') - x(\tau')
&\leq
y(\tau)- x(\tau)  + C_1 (\tau'- \tau)\,
2^{-\nu}
\\
\noalign{\medskip}
&\leq
y(\tau)- x(\tau)  + C_1
\displaystyle{\frac{x(\tau') - x(\tau)}{\lambda^{\min}} }
\cdot\big( w_k(\vfi_\nu(y(\tau')))- w_k(\vfi_\nu(x(\tau'))) \big)
\\
\noalign{\medskip}
&\leq
y(\tau)- x(\tau)  + C_1
\displaystyle{\frac{5 \rho'}{\lambda^{\min}} }
\cdot \big(y(\tau') - x(\tau')\big)
\end{array}
}
$$
which, because of (\ref{oleyncond}), implies
$$
y(\tau)- x(\tau) \geq
\bigg(
1-C_1\displaystyle{\frac{5 \rho'}{\lambda^{\min}}}
\bigg)\cdot \big(y(\tau') - x(\tau')\big)>0\,,
$$
proving (\ref{distraref2}).
\hfill$\Box$
\vsp

We next derive a regularity property enjoyed by 
general BV solutions
of Temple systems   
defined as limit of front tracking approximations,
which allows us to establish Proposition~\ref{theorem26}.
This is an extension of the regularity results
obtained in~\cite{dip, liu, brfiocco} for solution
with small total variation of general 
systems.
The arguments of the proof are quite similar
as for the corresponding result in~\cite{brfiocco},
but we will repeat some of them for completeness,
referring to~\cite{brfiocco}  
(see also \cite[Theorem~10.4]{br4})
for further details.

\begin{lemma}
\label{regularity}
Let  (\ref{311}) be a system of Temple class
satisfying the same assumptions as \linebreak in
Theorem \ref{theorem31}.
Consider a sequence
$u^\nu  
(t, \cdot) :[c,\, d]\mapsto \Gamma^\nu$, $t\in  [r,\,s]$,
of front tracking
approximate solutions  of the mixed problem for (\ref{311})
(constructed by the algorithm
of Section~\ref{subsec:existence}),
that converges  in $\lu$, as $\nu \to \infty$, to some function
$u(t, \cdot) : [c,\, d] \mapsto \Gamma$,
for every $t\in [r,\, s]\subset \R^+$.
Assume that
\be 
\tv (u^\nu (t, \cdot))\le M,
\qquad\quad
\tv (u^\nu (\cdot, x))\le M
\qquad\ \forall~t,x,\,\nu\,,
\label{br1}
\ee
for some constant $M>0$.
Then,
%
%
there exist a countable set of interaction points \linebreak
$\Theta\doteq \big\{ (\tau_l,\, x_l ); ~l\in
\NN\big\}\subset D\doteq [r,\,s]\times [c,\, d]$,
and a countable family of Lipschitz continuous shock curves
$\Upsilon \doteq \big\{ x=y_m(t);~t\in \,]r_m,\,s_m[,\>m\in
\NN\big\}$, such that the following hold.

\begin{itemize}
\item[(i)]
For each $m\in\NN$, and
for any $\tau\in \, ]r_m,\,s_m[$ with $(\tau, y_m(\tau))\not\in
\Theta$, there exist
the left and right  limits (\ref{limiti})
of $u$ at $(\tau, y_m(\tau))$
and the shock speed
$\dot y_m(\tau)$. Moreover, these limits
satisfy the Rankine Hugoniot relations (\ref{RH})
and the  Lax entropy inequality  (\ref{LAX}),
for some $i\in \{1,...,n\}$.

\item[(ii)]The map $u$ is continuous outside the set $\Theta \cup \Upsilon.$
\end{itemize}
\end{lemma}
\n\textsc{Proof.}

\n
{\bf 1.}
To establish $(i)$ we need to recall some technical tools
introduced in \cite{brfiocco} (see also  \cite[Theorem~10.4]{br4}).
For every front tracking solution $u^\nu$, we define  the {\it
interaction and cancellation measure} $\mu_\nu^{IC}$
that is a positive,  purely atomic measure on $D$,
concentrated on the set of points~$P$ where two or more
wave-fronts of $u^\nu$ interact. Namely,
 if the incoming
fronts at $P$ have size $\sigma_1, \dots , \sigma_\ell$
(w.r.t. the Riemann coordinates),
and belong to the families $i_1, \dots , i_\ell$ respectively,
we set
\be\mu_\nu^{IC}(P) \doteq
\sum_{\alpha, \beta} \big|\sigma_\alpha\,\sigma_\beta\big|+
\sum_{i}
\Bigg(
\sum_{\{i_\alpha\,;\,i_\alpha=i\}} |\sigma_\alpha|-
\bigg|\sum _{\{i_\alpha\,;\,i_\alpha=i\}} \sigma_\alpha\bigg|
\,\Bigg).
\label{misura}
\ee
Since $\mu_\nu^{IC}$ have a uniformly bounded total mass,
by possibly taking a subsequence we can  assume the
weak convergence
\be\mu_\nu^{IC} \rightharpoonup\mu^{IC} \label{convmeas}\ee
for some positive, purely atomic measure $\mu^{IC}$ on $D$.
Call $\Theta$ the countable set of atoms of~$\mu^{IC}$, i.e. set
$$\Theta \doteq \{P\in D ;~ \mu^{IC}(P)>0\}.$$
%
For every  approximate solution
$u^\nu$ taking values in $\Gamma^\nu$, $\nu \geq 1$,
and for any fixed $\eps \geq 2^{-\nu}$, 
by an {\it $\ve-$shock front of
the $i-$th family} in $u^\nu$ we mean
a polygonal line in $D$, with nodes
$(\tau_0,x_0),...,(\tau_N,x_N),$ having the following properties.
\begin{itemize}
\item[\textsc{(I)}] The nodes $(\tau_h, x_h)$
are interaction points or lie on the boundary of $D$,
and the sequence of times is increasing $\tau_0< \tau_1< \cdots<\tau_N$\,.

\item[\textsc{(II)}] Along each segment joining
$(\tau_{h-1},x_{h-1})$ with $(\tau_{h},x_{h})$, the function  $u^\nu$ has an $i-$shock
with strength $|\sigma_{h}|\ge \eps$.

\item[\textsc{(III)}]
For $h<N$,
if two (or more) 
incoming $i-$shocks of strength
$\ge \ve$ interact at the node~$(\tau_{h},x_{h})$, then the
shock coming from $(\tau_{h-1},x_{h-1})$ has the larger speed, i.e. is the one
coming from the left.

\end{itemize}

\n An $\ve-$shock front which is maximal with respect to the set theoretical
inclusion will be called a {\it maximal $\ve-$shock front}.
Observe that, because of \textsc{(III)}, two maximal $\ve-$shock fronts of the same family
either are disjoint or coincide.
Moreover, by (\ref{br1}), the number of maximal $\ve-$shock front
that starts at the boundary of $D$ is uniformly bounded by $3M/\eps$.
On the other hand, the special geometric features of Temple class
systems guarantee that no new shock front can arise in the interior
of $D$. Indeed, the coinciding shock and rarefaction assumption
together with 
the existence of
Riemann invariants prevents the creation of shocks 
of other families
than the ones of the incoming fronts at any interaction point. 
Therefore, for fixed $\eps>0$, and  $i\in\{1,\dots ,n\}$,
the number of maximal $\ve-$shock front of the $i$-th family
remains uniformly bounded by $M_\eps\doteq 3M/\eps$ in all $u^\nu, \, \nu\geq 1$.
Denote such curves by
$$
y_{\nu,m}^\eps:[t^{\eps,-}_{\nu,m},\,t^{\eps,+}_{\nu,m}]
\mapsto \R, \qquad m=1,..., M_\eps\,.
$$ 
By possibly extracting a further subsequence,
we can assume the convergence
$$
y_{\nu,m}^\eps(\cdot)\longrightarrow y^\eps_{m}(\cdot), \qquad\quad
t^{\eps,\pm}_{\nu,m}\longrightarrow t^{\eps,\pm}_{m},\
\qquad m=1,..., M_\eps\,,
$$
for  some Lipschitz continuous paths 
$y_{m}^\eps:[t^{\eps,-}_{m},\,t^{\eps,+}_{m} ]\mapsto \R, \>
m=1,..., M_\eps.$
Repeating this construction in connection with a 
sequence $\eps_k \to 0$, and taking the union of all the paths
thus obtained, we find, for each 
characteristic family $i\in\{1,\dots,n\}$,
a countable family of Lipschitz continuous  curves 
$y_{m}:[t^-_{m},t^+_{m} ]\mapsto \R, \>
m\in\NN$. Call $\Upsilon$ the union of all such curves.

%
\vsp

\n
{\bf 2.}
Consider now 
a point $P=(\tau,y_m(\tau))\not\in\Theta$ along a curve
$y_m\in \Upsilon$ of a family $i\in\{1,\dots,n\}$. 
Notice that, by construction, and because of (\ref{convmeas}), 
no curve in $\Upsilon$ can cross $y_m$ at $P$.
Moreover, 
by (\ref{br1}), the function $u(\tau, \cdot)$ has bounded variation,
and hence there exist the limits 
\be
\lim\limits_{x\rightarrow  y_m(\tau)-} u(\tau,x)
\doteq u^-,\qquad
\lim\limits_{x\rightarrow  y_m(\tau)+} u(\tau,x)
\doteq u^+.
\label{limiti2}
\ee
We claim that also the limits (\ref{limiti}) exist,
and thus coincide with those in (\ref{limiti2}).
To this end observe that, by construction, 
there exist a sequence of shocks curves $y_{\nu,m}$
of the $i$-th family
converging to $y_{m}$, along which each approximate solution $u^\nu$
has  a jump of strength~$\ge \ve^*$, for some $\eps^*>0$.
Then, relying on the assumption 
\be
\mu^{IC}(\{P\})=0\,,
\label{zeromass}
\ee
and letting $B(P,r)$ denote the ball centerd at $P$ with radius $r$,
one can establish the limits
\be
\lim\limits_{r\rightarrow 0+} 
\limsup\limits_{\nu \rightarrow +\infty}
\Bigg(\sup_{(t,x)\in B(P,r),
\>
x<y_{\nu,m}(t)} \big| u^\nu(t,x)-u^-\big|
\Bigg)=0, \label{lims1}
\ee
\be
\lim\limits_{r\rightarrow 0+} 
\limsup\limits_{\nu \rightarrow +\infty}
\Bigg(\sup_{
(t,x)\in B(P,r)
\>
x>y_{\nu,m}(t)
} \big| u^\nu(t,x)-u^-\big|
\Bigg)=0, \label{lims2}
\ee
which clearly yield (\ref{limiti}). Indeed, if for example (\ref{lims1})
do not hold, by possibly taking  a subsequence
we would find $\eps>0$ and points $P_\nu\doteq (t_\nu,\,\xi_\nu) \to P$
on the left of $y_{\nu,m}$ such that
$$\big|u^\nu(t_\nu,\,\xi_\nu)-u^-\big|\geq \eps
\qquad\quad \forall~\nu.$$
On the other hand, by the first limit in (\ref{limiti2}),
and since $u^\nu(\tau ,x) \to u(\tau,x)$ for a.e. $x\in [\alpha, \beta]$,
we could also find points $Q_\nu\doteq (\tau, \, \xi_\nu')\to P$ on the left
of $y_{\nu,m}$ such that
$$
u^\nu(\tau, \, \xi_\nu') \to u^-,
\qquad\qquad
\frac{|\xi_\nu-\xi_\nu'|}{|t_\nu-\tau|}> \lambda^{\max}
\qquad \forall~\nu\,,
$$
where $\lambda^{\max}$ denotes the maximum speed at (\ref{uppbound}).
But then, for each solution $u^\nu$, the segment
$\overline{P_\nu\,Q_\nu}$ would be crossed
by an amount of waves of strength $\geq \eps$.
Hence, by strict hyperbolicity and genuine nonlinearity,
this would generate a uniformly positive amount 
of interaction and cancellation within an arbitrary 
small neighborhood of $P$ 
(see. \cite[Theorem 10.4-Step 5]{br4})
which, by the definition (\ref{misura}),
and because of (\ref{convmeas}),
contradicts the assumption (\ref{zeromass}).

\vsp

To complete the proof of (i) 
observe that, by construction,
the states $u^-_{\nu,m}(\tau), \, u^+_{\nu,m}(\tau)$ to the left
and to the right of the jump in $u^\nu$ at $y_{\nu,m}(\tau)$
satisfy the Rankine  Hugoniot conditions. Thus,
 relying on (\ref{lims1})-(\ref{lims2}), 
and on the convergence $y_{\nu,m} \to y_\nu$, one 
deduces (\ref{RH}).
The proof of (ii) can be established with the same type
of arguments (cfr. \cite[Theorem 10.4-Step 8]{br4}).
\hfill$\Box$
\vsp

As an immediate consequence of Lemma~\ref{regularity}, we derive
Proposition~\ref{theorem26} stated in Section~\ref{subsec:existence}.
\vsp

\n
\textsc{Proof of Proposition~\ref{theorem26}.}
Consider a sequence 
$u^\nu(t, \cdot) : [a, b] \mapsto \Gamma^\nu$
of front tracking approximate solutions
of the mixed problem for~(\ref{311})
on the region $\Omega_T\doteq [0,T]\times[a,b]$,
that converges  in~$\lu$, as $\nu \to \infty$, to some function
$u(t, \cdot) : [a, b] \mapsto \Gamma$,
for every $t\in [0,T]$.
Observe that, by Theorem~\ref{theorem31}
 one can find another sequence $\{v^\nu\}_{\nu \geq 1}$
of approximate solutions of (\ref{311}) on the region $\Omega_T$, whose initial
and boundary data have a number of shocks $N_\nu \leq \nu$
for each characteristic family, and such that 
$$
\big\|u^\nu(t, \cdot)-v^\nu(t, \cdot)\big\|_{\lu([a,b])}\leq 1/\nu
\qquad\quad\forall~t\in[1/\nu,\,T]\,.
$$
Then,
thanks to the Oleinik estimates (\ref{oleinappr1})-(\ref{oleinappr2}),
and because all~$v^\nu$
take values in the compact set~(\ref{3212}),
there will be,
for every fixed $\eps>0$, some constant~$M_\eps>0$
such that
\be
{\setlength\arraycolsep{2pt}
\begin{array}{rr}
\tv \big\{v^\nu(t,\, \cdot)~; ~[a+\eps,\, b-\eps]\big\} \leq& M_\eps
\qquad\quad \forall~t\in [\eps,\, T]\,,
\hskip 32pt
\\
\noalign{\smallskip}
\tv \big\{v^\nu(\cdot,\, x)~; ~[\eps,\, T]\big\} \leq& M_\eps
 \qquad\quad \forall~x \in [a+\eps,\, b-\eps] \,,
\end{array}
}
\qquad\ \forall~\nu\in\NN.
\label{tvbound}
\ee
%
%
%
Thus, writing $\Omega_T$
as the countable union
$$
\Omega_T=\cup_k D_k,\qquad\quad D_k\doteq  [1/k,\, T] \times  [a+(1/k),\, b-(1/k)],
$$
and applying Lemma~\ref{regularity} to each
sequence of maps $v^\nu_k\doteq v^\nu\!\restriction_{D_k},\, \nu\geq 1,$
defined as the restriction of $v^\nu$ to the domain $D_k$, we
clearly reach the conclusion of Proposition~\ref{theorem26}.
\hfill$\Box$
\vsp

We are now in the position to establish Theorem~\ref{theorem25},
relying on Proposition~\ref{theorem26} and on~Theorem~\ref{theorem32}.
\vsp

\n
\textsc{Proof of Theorem~\ref{theorem25}.}
Let
$u^\nu(t, \cdot) : [a, b] \mapsto \Gamma^\nu$
be a sequence
of front tracking approximate solutions
of the mixed problem for~(\ref{311})
on the region $\Omega_T\doteq [0,T]\times[a,b]$,
that converges  in $\lu$, as $\nu \to \infty$, to some function
$u(t, \cdot) : [a, b] \mapsto \Gamma$,
for every $t\in [0,T]$.
Since, by construction, each $u^\nu$ is a weak solution
of~(\ref{311}), and because
$u^\nu(0, \cdot) \to u(0, \cdot)=\overline u$,
also the limit function $u$ is a weak solution of
the Cauchy problem (\ref{311})-(\ref{312}) on the region~$\Omega_T$.
Moreover, applying Proposition~\ref{theorem26},
we deduce that $u$ admits at $t=0$ and at $x=a, \, x=b$ 
the left and right limits
 stated in~Theorem~\ref{theorem25}.
On the other hand, by the same arguments
used in the proof of Proposition~\ref{theorem26},
we may assume that the  initial
and boundary data of each approximate 
solution~$u^\nu$ have at most $N_\nu \leq \nu$
shocks for every characteristic family.
Then, letting $\nu \to \infty$ in (\ref{oleinappr1})-(\ref{oleinappr2}),
by the lower semicontinuity
of the total variation we find that $u$
satisfies  the entropy conditions (\ref{330})-(\ref{328b})
on the decay of positive waves.
It follows that $u$ is an
entropy weak solution of the mixed problem 
(\ref{311})-(\ref{314}) according with 
Definition~\ref{definition32}.
Hence, observing that by constructin the map
$(t,x) \to \big(u(t, \cdot),\,u(\cdot, x)\big)$
takes values within the 
domain $\D_T$ defined in (\ref{domainT}),
and applying Theorem~\ref{theorem32}, we
deduce that (\ref{3218}) is verified. 
\hfill$\Box$
\vsp

\section{Proof of Theorems \ref{theorem1}-\ref{theorem2}}
\label{section34}
\setcounter{equation}{0}
\indent

\n
\textsc{Proof of Theorem~\ref{theorem1}.}  \
%
 We shall first prove that, for every fixed 
 $\overline\tau>0$, there exists some 
constant $\rho=\rho(\overline\tau)>0$
so that (\ref{t2}) holds.
Given $\widetilde u_a\in\U^{\infty}_{\tau}, \,
\widetilde u_b\in\U^{\infty}_{\tau}$, \, 
$\tau \geq \overline\tau$, 
let $u=u(t,x)$ be an entropy weak
solution of (\ref{311})-(\ref{314}) on
the region $[0, \tau]\times [a, b]$ according with
Definition~\ref{definition32}. Then, the Oleinik-type estimates
(\ref{331}) on the decay of positive waves
imply that, for $i\in \{p+1,...,n\}$,
 $\tau\geq\overline\tau$, and for  a.e. $a <x<y<b$,  
there holds
\be
{\setlength\arraycolsep{2pt}
\begin{array}{ll}
\displaystyle{
\frac{w_i(\tau,y)-w_i(\tau,x)}{y-x}}
&\leq 
C\cdot 
\displaystyle{\bigg\{
\frac{y-x}{\tau}+\log\bigg(\frac{y-a}{x-a}\bigg)\bigg\}}
\\
\noalign{\medskip}
&\le (b-a)\,C\cdot 
\displaystyle{\bigg\{
\frac{1}{\overline\tau}+\frac{1}{x-a}\bigg\}}
\\
\noalign{\medskip}
&\leq 
\displaystyle{
\frac{C\,(b-a)\big((b-a)+\overline\tau\big)}{\overline\tau}
\cdot\frac{1}{x-a}}.
\end{array}}
\label{31}
\ee
Clearly, with the same computations, relying on 
the Oleinik-type estimates
(\ref{330}), we deduce
that, for $i\in \{1,...,p\}$,
 $\tau\geq\overline\tau$, and for  a.e. $a <x<y<b$,  
there holds
\be
\displaystyle{
\frac{w_i(\tau,y)-w_i(\tau,x)}{y-x}}
\le
\displaystyle{
\frac{C\,(b-a)\big((b-a)+\overline\tau\big)}{\overline\tau}
\cdot\frac{1}{b-y}}.
\label{32}
\ee
Hence, taking 
\be
\rho
\ge  \displaystyle{
\frac{C\,(b-a)\big((b-a)+\overline\tau\big)}{\overline\tau}}
\ee
from (\ref{31})-(\ref{32}) we derive 
$u (\tau, \cdot) \in K^{\rho}$,  
which proves $(\ref{t2})$. 
\vsp

%

Concerning the second statement of the theorem, we
will show that, letting $\lambda^{\min}$, $\rho'$, be
the minimum speed in (\ref{uppbound}),
and the constant (\ref{oleyncond}) of Lemma~\ref{lemma1}, and taking
\be
T\doteq
\displaystyle{\frac{4\,(b-a)}{\lambda^{\min}}}
\label{time}
\ee
%
the relation (\ref{t1})
is verified,
i.e. that, given $\vfi\in K^{\rho'}$, 
and~$\tau>T$, there exist $\widetilde u_a\in\U^{\infty}_{\tau}$, \,
$\widetilde u_b\in\U^{\infty}_{\tau}$, and a solution
$u(t,x)$ of 
(\ref{311})-(\ref{314}) on $[0,\tau]\times[a,b]$
(according with Definition~\ref{definition32}),
such that $u(\tau, \cdot) \equiv \vfi.$
Notice that, 
by Remark~\ref{remark21}, we may assume that 
$w_i(\vfi(x))$ admits left and right limits
in any point $x\in ]a, \, b[$,
and that $w_i(\vfi(x))=w_i(\vfi(x^+))\doteq 
\lim_{\xi\to x^+} w_i(\vfi(\xi))$, for
all $i \in \{1,\dots,n\}$.
The proof is devided in two steps.
\vsp

\n{\bf Step 1. Backward
construction of front tracking approximations.} \
Letting $\rho'>0$ be the constant  in (\ref{oleyncond}),
consider a sequence $\{\vfi_\nu\}_{\nu\geq 1}$  of
(right continuous) piecewise constant maps in~$K_\nu^{\rho'}$, 
satisfying the conditions a)-b) of Lemma~\ref{lemma2},
and take a piecewise constant approximation
$\overline u^\nu:[a,b]\mapsto \Gamma^{\nu}$ 
of the initial data $\overline u$, so that
$\overline u^\nu \to \overline u$ in $\lu$.
Given $\tau>T$ ($T$ being the time defined in (\ref{time})),
for each $\nu \geq 1$, we will construct here a
front tracking solution $u^\nu(t,x)$ of~(\ref{311}) 
on the region $[0,\tau]\times [a, b]$,
with initial data $u^\nu(0,\cdot)=\overline u^\nu$,
so that
\be
u^\nu(\tau, \cdot\,)= \vfi_\nu\,.
\label{finalcond}
\ee
%
This goal is accomplished by proving the following two lemmas.
%
\begin{lemma} \label{lemma4}
Let $T,\, \rho'>0$ be the constants  in (\ref{time})
and (\ref{oleyncond}).
Then, for  every (right continuous) 
$\vfi_\nu\in K_\nu^{\rho'}$,  $\nu \geq 1$,
satisfying the condition a) of Lemma~\ref{lemma2},
and for any $\tau>T$, there exists a front tracking solution $u^\nu(t,x)$
of~(\ref{311}) 
on the region $[(3/4)T,\, \tau]\times [a, b]$,
with boundary  data
$\widetilde u^\nu_a \doteq u^\nu(\cdot,\,a), \,
\widetilde u^\nu_b \doteq u^\nu(\cdot,\,b)
\in \li([(3/4)T,\,\tau],\,\Gamma^\nu)$, so that
\be
u^\nu\big((3/4)T,\, x\big) \equiv \omega,
\qquad\quad
u^\nu\big(\tau,\, x\big) =\vfi_\nu(x)\,,
\qquad\quad\forall~x\in [a,b]\,,
\label{4.7}
\ee
for some constant state $\omega\in \Gamma^\nu$.
\end{lemma}

\n
\textsc{Proof.}
Given $\tau>T$, and
$\vfi_\nu\in K_\nu^{\rho'}$, $\nu \geq 1$,
satisfying the condition a) of Lemma~\ref{lemma2},
we will use the algorithm described in Section~\ref{subsec:existence}
to construct backward in time a front tracking solution
that takes value $\vfi_\nu$ at time $\tau$.
To this end, we first observe that
according with the algorithm of Section~\ref{subsec:existence}
we can always 
construct the
backward solution of a Riemann problem
with terminal data
\be
u(t,x)=
\left\{
\!\begin{array}{ll}
u^L\quad &\hbox{if}\qquad x<\xi\,,
\\
u^R\quad &\hbox{if}\qquad x>\xi\,,
\end{array}
\right.
\label{tdrp}
\ee
if the  the terminal states 
$u^L, \, u^R\in \Gamma^\nu$
have Riemann coordinates 
$$
w(u^L)\doteq w^L=(w^L_1,\ldots,w_n^L),\qquad 
w(u^R)\doteq w^R=(w^R_1,\ldots,w_n^R)
$$
that satisfy
\be
w_i^L < w_i^R
\qquad \Longrightarrow
\qquad
w_i^R = w_i^L + 2^{-\nu}
\qquad\quad \forall~i\,.
\label{rpcond}
\ee
Indeed, if we consider the intermediate states
\be
z^i=
\left\{
\!\begin{array}{ll}
u^L\quad &\hbox{if}\qquad i=0\,,
\\
u(w^L_1,\ldots,w^L_{n-i}\,,\,w^R_{n-i+1},\ldots,w^R_n)\quad 
&\hbox{if}\qquad 0<i<n\,,
\\
u^R\quad &\hbox{if}\qquad i=n\,,
\end{array}
\right.
\ee
we realize that,
because of (\ref{rpcond}), the solution
of every 
Riemann problem with initial \linebreak states~$(z^{i-1}, z^i)$ (defined as
in Section~\ref{subsec:existence})
contains only a single front.
Thus,
we can construct
the 
solution to the Riemann problem 
with terminal data~(\ref{tdrp})
in a backward \linebreak neighborhood of~$(t,\xi)$
by piecing together the  solutions to the simple
Riemann problems $(z^{i-1}, z^i)$, $i=1,\dots,n$.

A front tracking solution $u^\nu$
can now be constructed
backward in time starting at~$t=\tau$,
and piecing together the backward solutions
of the Riemann problems determined by
the jumps in~$\vfi_\nu$.
The resulting piecewise constant function $u^\nu(\tau-, \,\cdot)$
is then
prolonged for~$t<\tau$ tracing backward the incoming fronts 
at $t=\tau$, up to the first time $\tau_1<\tau$ 
at which two or more discontinuities cross
in the interior of $\Omega$.
Observe that,
since $u^\nu$ is a front tracking solution
constructed by the algorithm of Section~\ref{subsec:existence}
on the region $[\tau_1,\tau]\times [a,b]$, we can apply
Lemma~\ref{lemma3}. Hence, 
it follows that the left and right states
of the jumps occuring in
$u^\nu(\tau_1,\,\cdot)$ satisfy condition (\ref{rpcond}),
because (\ref{distraref2})
guarantees that two (or more) adjacent rarefaction fronts of the same family
cannot cross at time $\tau_1$.
We then solve backward the resulting Riemann problems
applying again the above procedure.
This determines the solution $u^\nu(t,\cdot)$ untill the
time $\tau_2<\tau_1$ at which another intersection between its
fronts takes place in the interior of $\Omega$,
and so on (see figure 1a).

\parbox{5cm}{
 \begin{center}
\leavevmode \epsfxsize=2.5in \epsfbox{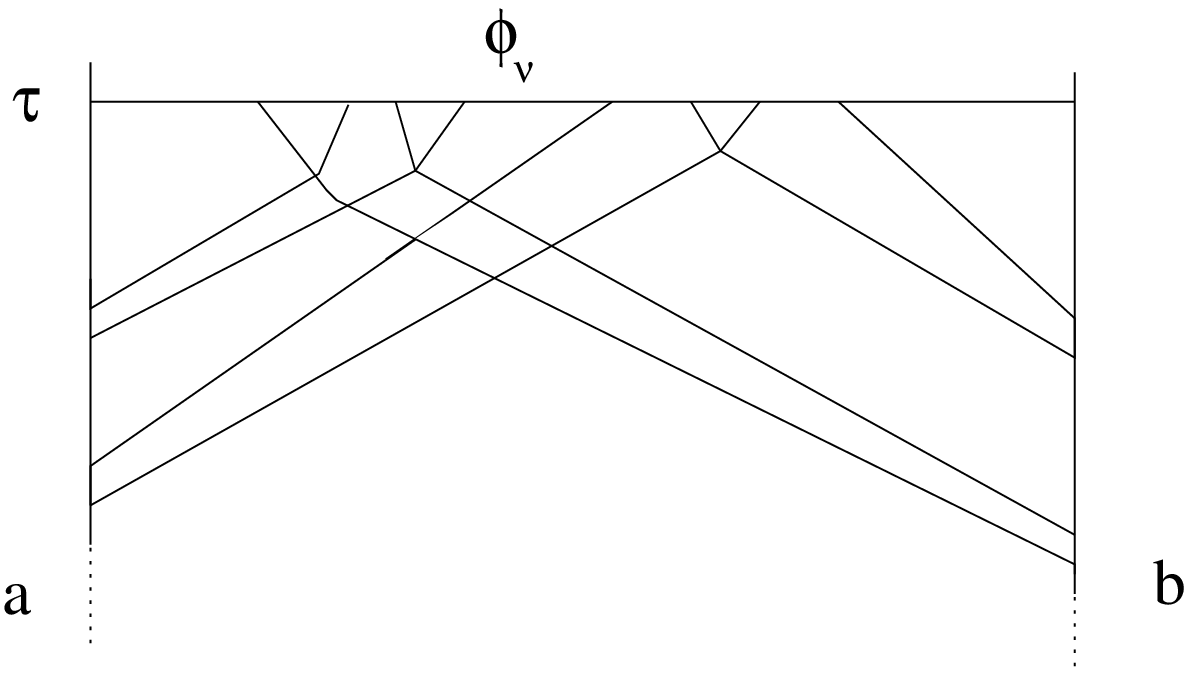}     
\end{center}     }
\hfill
\parbox{5cm}{  
 \begin{center}
\leavevmode \epsfxsize=1.7in \epsfbox{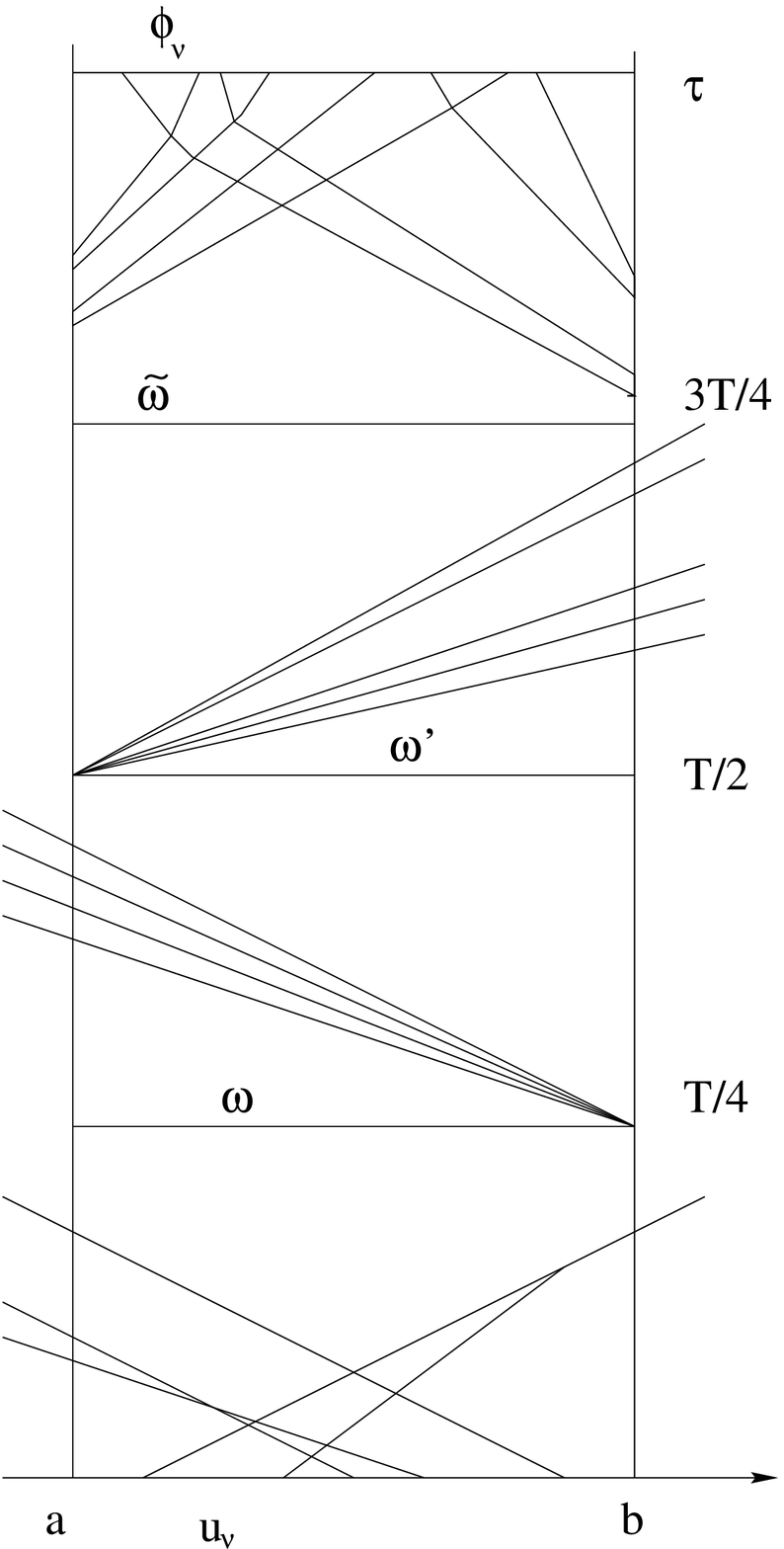}
\end{center}         }
\begin{center} 
~~~~~figure 1a~~~~~~~~
~~~~~~~~~~~~~~~~~~~~~~~~~~~~~~~~~~~~~~~~~~~~~~~~~~figure 1b
   \end{center}

With this construction we define a front tracking solution 
$u^\nu(t,x)$ on the whole \linebreak region $[(3/4)T,\, \tau]\times [a, b]$,
that verifies the first equality in (\ref{4.7}),
and corresponds to the boundary data
$\widetilde u^\nu_a\doteq u^\nu(\cdot,a), \,
\widetilde u^\nu_b \doteq u^\nu(\cdot,b)\in
\li([(3/4)T,\,\tau],\,\Gamma^\nu)$.
Clearly, the total number of wave-fronts in $u^\nu(t, \, \cdot)$
decreases, as $t \downarrow (3/4)T$, whenever a (backward) front crosses
the boundary points $x=a, \, x=b$.
Since (\ref{uppbound}) implies that 
the maximum time taken by fronts of $u^\nu$ to
cross the interval~$[a,b]$
is $(b-a)/\lambda^{\min}$, the definition~(\ref{time}) of $T$
guarantees that all the (backward) fronts of $u^\nu$ will hit 
the boundaries $x=a, \, x=b$ within 
some time $\tau'\in\,](3/4)T,\,\tau[$\,,
which shows that also the second equality in (\ref{4.7})
is verified, thus completing the proof.
\hfill$\Box$
\vskip 5pt

\begin{lemma} \label{lemma5}
Let $T>0$ be the constant in (\ref{time}).
Then, for any piecewise constant function \linebreak 
$\overline u^\nu\in \li([a,b],\Gamma^\nu)$,
and for every state $\omega\in \Gamma^\nu$,
there exists a
front tracking solution \linebreak $u^\nu(t,x)$ of~(\ref{311}) 
on the region $[0,\,(3/4)\,T]\times [a, b]$,
corresponding to some boundary  data \linebreak 
$\widetilde u^\nu_a, \,
\widetilde u^\nu_b \in \li([0,\,(3/4)\,T],\Gamma^\nu)$,
so that
\be
u^\nu\big(0,\, x\big) =\overline u^\nu(x),
%
%
\qquad \quad  
u^\nu\big((3/4)\,T,\, x\big) \equiv\omega\,,
\qquad\quad\forall~x\in [a,b]\,.
\label{finalcondb}
\ee
%
%
\end{lemma}

\n
\textsc{Proof.}
The approximate solution $u^\nu$ is constructed as follows.
By Remark~\ref{remark24}, \linebreak for $t\in [0, T/4]$, we can define
$u^\nu(t,x)$ as the restriction to the region $[0, T/4]\times [a,b]$
of the front tracking solution 
to the Cauchy problem for (\ref{311}), with initial data
$$\overline u(x) = \left\{
\!\!\begin{array}{ll}
\ov u^\nu(a+) & \>\textrm{if \ \ $ x<a$},\\ 
\ov u^\nu(x) & \>\textrm{if \ \ $a \le x \le b$},\\
\ov u^\nu(b-) & \>\textrm{if \ \ $ x > b$},
\end{array} \right.
$$
(constructed as in \cite{brgoa1} with the same type of
algorithm described in Section~\ref{subsec:existence}).
Observe that, since $u^\nu$ contains only fronts originated
at the points of the segment $\{(0,x)\,;\,x\in[a,b]\}$,
because of (\ref{uppbound}), (\ref{time}) these wave-fronts 
cross the whole interval $[a, b]$ and exit from the
boundaries $x=a, \, x=b$ before time $T/4$ (see figure 1b).
Hence, there will be some state $\omega'\in \Gamma^\nu$
such that
\be
u^\nu\big(T/4,\, x\big) \equiv\omega'
\qquad\ \ \forall~x\in [a,b]\,.
\label{intermstates}
\ee
Thus, introducing the intermediate state 
$$
\widetilde \omega\doteq \big(\omega_1, \dots, \omega_p, \,
\omega_{p+1}', \dots , \omega_n' \big)
$$
between $\omega'$ and $\omega$, 
we will define $u^\nu(t,x)$, for $t\in [T/4, \, T/2]$, as the restriction
to the region $[T/4, T/2]\times [a,b]$
of the approximate solution to the
Riemann problem for (\ref{311}), with initial
data
\be
u^\nu(T/4,x)=
\left\{
\!\begin{array}{ll}
u(\omega')\quad &\hbox{if}\qquad x<b\,,
\\
u(\widetilde \omega)\quad &\hbox{if}\qquad x>b\,,
\end{array}
\right.
\label{rp2}
\ee
while, for $t\in [T/2, \, (3/4)T]$, we will let 
$u^\nu(t,x)$ be 
the restriction
to the region $[T/2, \, (3/4)T]\times [a,b]$
of the approximate solution to
the Riemann problem for (\ref{311}), with initial
data
\be
u^\nu(T/2,x)=
\left\{
\!\begin{array}{ll}
u(\omega)\quad &\hbox{if}\qquad x<a\,,
\\
u(\widetilde \omega)\quad &\hbox{if}\qquad x>a\,.
\end{array}
\right.
\label{rp3}
\ee
By the definition of $\widetilde \omega$,
and because of (\ref{uppbound}), (\ref{time}),
on $[T/4, \, T/2]$ the solution of the Riemann problems
with initial data (\ref{rp2}) contains 
only wave-fronts originated
at the point~$(T/4, \, b)$, 
that cross the whole interval $[a, b]$ and exit from the
boundary $x=a$ before time~$T/2$.
Similarly, still by (\ref{uppbound}), (\ref{time}),
for $t\in  [T/2, \, (3/4)T]$
the solution of the Riemann problem
with initial data (\ref{rp3}), 
contains only wave-fronts originated
at $(T/2, a)$, that cross the whole interval~$[a, b]$, 
and exit from the
boundary $x=b$ before time~$(3/4)T$ (see figure 1b).
Hence, $u^\nu(t,x)$
is a front-tracking solution 
defined on the whole region 
$[0,\,(3/4)\,T]\times [a, b]$,
 that corresponds to the boundary data
$\widetilde u^\nu_a \doteq u^\nu(\cdot,a)$, \,
$\widetilde u^\nu_b \doteq u^\nu(\cdot,b)
\in \li([0,\,(3/4)\,T],\Gamma^\nu)$,
and verifies 
the conditions~(\ref{finalcondb}).
\hfill$\Box$

\vsp

\n{\bf Step 2. Convergence of the approximate solutions.}
By Step~1, for a given $\vfi\in K^{\rho'}$
(with~$\rho'$ as in (\ref{oleyncond})),
 we have found a sequence of initial data $\overline u^\nu$,
and of boundary data 
$\widetilde u^\nu_a, \, \widetilde u^\nu_b\in \U^\infty_{\tau}$,
so that, letting $u^\nu(\tau, \cdot)\doteq E_\tau^\nu
(\overline u^\nu, \, \widetilde u^\nu_a, \, \widetilde u^\nu_b)$
be the corresponding front tracking solution, there holds
\be 
\overline u^\nu \to \overline u,
\qquad\quad 
u^\nu(\tau, \cdot)\rightarrow \vfi
\qquad\hbox{in}\qquad \lu([a, b]).
\label{approxconv2}
\ee
By the same arguments
used in the proof of Proposition~\ref{theorem26},
we may assume that the  initial
and boundary data of each approximate 
solution~$u^\nu$ have at most $N_\nu \leq \nu$
shocks for every characteristic family.
Then, thanks to the Oleinik-type estimates
(\ref{oleinappr1}),
and because $u^\nu$ are uniformly bounded since they
take values in the compact set~(\ref{3212}),
for every fixed $\eps>0$, 
there will be some constant~$C_\eps>0$
such that
\be
{\setlength\arraycolsep{2pt}
\begin{array}{rr}
\tv\big\{u^\nu (t,\cdot\,)~;~[a+\eps,\,b-\eps]\big\}\leq& C_\eps
\qquad\forall~t\in [\eps,\tau]\,,\hskip 35pt
\\
\noalign{\medskip}
\displaystyle{
\int_{a+\eps}^{b-\eps}|u^\nu(t,x)-u^\nu(s,x)|~dx}\leq& C_\eps |t-s|
\qquad \forall~t,s\in [\eps,\tau]\,,
\end{array}
}
\qquad\ \forall~\nu\in \NN\,.
\label{tcont}
\ee
Hence, applying  Helly's Theorem,  we deduce that
there exists a subsequence $\{u^{\nu_j}\}_{j\geq 0}$ that
converges in $\lu([a,b],\,\Gamma)$
to some function $u_{\eps}(t,\cdot),$
for any $t\in [\eps,\,\tau].$
Therefore, repeating the same construction
in connection with a sequence
$\eps_k \to 0+,$
and using a diagonal procedure, we obtain a subsequence
$\{u^{\nu'}(t,\cdot\,)\}_{\nu'\geq 0}$ that converges
in $\lu([a,b], \Gamma)$
to some function $u(t,\cdot),$ 
for any $t\in [0,\,\tau].$
Then, by Theorem~\ref{theorem25}, 
there holds (\ref{3218}),
with 
$\widetilde u_a \doteq
u(\cdot,a), \, \widetilde u_b\doteq u(\cdot,b)\in \U_\tau^\infty$, 
while~(\ref{approxconv2}) implies $u(\tau, \cdot)= \vfi$,
which shows~$\vfi \in \A(\tau)$. This completes the proof 
of Theorem~\ref{theorem1}.
\hfill$\Box$
\vsp

We next establish the compactness of the attainable set (\ref{attset2})
stated in Theorem~\ref{theorem2}.
The proof is quite similar to that of \cite[Theorem~2.3]{A-G}.
We repeat it for completeness.
\vsp

\n\textsc{Proof of Theorem~\ref{theorem2}.} 
Fix $T>0$, and  
consider a sequence $\{u^\nu\}_{\nu\geq 0}$ of entropy weak solutions 
to the mixed problem for (\ref{311}) 
on $\Omega_T\doteq [0, T]\times [a,b]$
(according with Definition~\ref{definition32}), with 
a fixed initial data $\overline u \in \li([a,b], \Gamma)$.
Since all $u^\nu$ are uniformly bounded, and because 
of the Oleinik-type estimates~(\ref{330})-(\ref{331}), one can find, 
for every $\eps>0$, some constant~$C_\eps>0$
so that (\ref{tcont}) holds.
Thus, with the same arguments used in {\bf Step 2} 
of the previous proof, we can construct a subsequence 
$\{u^{\nu'}\}_{\nu'\geq 0}$ so
that, for any $t\in [0,\,T]$, $u^{\nu'}(t,\cdot\,)$ converges
in $\lu$
to some function~$u(t,\cdot),$ 
 which is continuous as a map
from $]0, T]$ into $\lu([a,b], \Gamma)$,
and 
satisfies  the entropy conditions (\ref{330})-(\ref{328b})
on the decay of positive waves.
On the other hand, the weak traces~$\Psi^{\nu'}_a, \, \Psi^{\nu'}_b$
of the fluxes
$f(u^{\nu'})$ at the boundaries $x=a, \, x=b$ 
are uniformly bounded, and hence are weak$^*$ relatively compact
in $\li ([0,T])$. 
Thus, by possibly taking a further subsequence, we have
\be
\Psi_a^{\nu'} \mathop{\rightharpoonup}^* \Psi_a\,, \qquad\quad\Psi_b^{\nu'}
\mathop{\rightharpoonup}^* \Psi_b\qquad\hbox{in }\quad \li([0,T]),
\label{weaktracelim} 
\ee
for some maps $\Psi_a,\> \Psi_b\, \in\li ([0,T])$.
Notice that, by the properties of the Riemann invariants,
the set $f(\Gamma)$ is closed and convex, and hence 
also the weak limits  $\Psi_a,\> \Psi_b$ 
take values in $f(\Gamma)$.
Moreover, since each $u^\nu$ is a distributional solution
of (\ref{311})-(\ref{312}) on  $\Omega_T$, also the limit function $u$
is a distributional solution of the Cauchy problem 
(\ref{311})-(\ref{312}) on the region $\Omega_T$.
Then, setting $\widetilde u_a \doteq f^{-1}\circ \Psi_a, \, 
\widetilde u_b \doteq f^{-1}\circ \Psi_b$,
it follows that $u$ is an
entropy weak solution of the mixed 
problem~(\ref{311})-(\ref{314}) (with boundary data in $\U^\infty_T$)
according with 
Definition~\ref{definition32},
which shows that $u(T, \cdot) \in \A(T)$. This completes the proof 
of Theorem~\ref{theorem2}.
\hfill$\Box$
\vsp
 
If we take in consideration only solutions to the mixed problem
(\ref{311})-(\ref{314}) that are trajectories
of the flow map $E$ obtained in Theorem~\ref{theorem31}
(which, in particular, admit a strong $\lu$ trace at the boundaries
$x=a, x=b$), we are lead to study
the set of attainable profiles
%
%
\be
\A_E(T)\doteq
\big\{
E_T  (\,\overline u,\,\widetilde u_a\,,\,\widetilde u_b\,)~; ~~
\quad \widetilde u_a,\, \widetilde u_b \in 
{\bf L}^{\infty}([0,T],\Gamma)
\big\}\,.
\label{attset3}
\ee
Since $\A_E(T)\subset \A(T)$, and by the proof of 
Theorem~\ref{theorem1}, it clearly follows that
the characterization of the set $\A(T)$ provided by
the  inclusions
(\ref{t2})-(\ref{t1}) 
of Theorem~\ref{theorem1} 
holds also for~$\A_E(T)$. 
Concerning the  compactness of the set $\A_E(T)$,
observe that, given any sequence of exact
solutions 
$u^\nu(t, \cdot) \doteq 
E_t (\,\overline u^\nu,\,\widetilde u_a^\nu\,,\,\widetilde u_b^\nu\,)$,
$\nu \geq 1$,
by Theorem~\ref{theorem31}
one can find another sequence of approximate solutions
$v^\nu(t, \cdot)$ constructed by the front tracking algorithm
of Section~\ref{subsec:existence}, so that
%
$$
\big\|u^\nu(t, \cdot)-v^\nu(t, \cdot)\big\|_{\lu([a,b])}\leq 1/\nu
\qquad\quad\forall~t\in[1/\nu,\,T]\,.
$$
Therefore, 
relying on
the regularity property of a solution obtained
as limit of front tracking approximations provided
by Theorem~\ref{theorem25},
with the same arguments
used in the proof of Theorem~\ref{theorem2} 
one can establish also the compactness of the set $\A_E(T)$.

\begin{theorem} \label{theorem3}
Under the same assumptions of Theorem 1, the set $\A_E(T)$ is 
a compact subset of~$\lu ([a,b],\,\Gamma)$ for
each $T>0.$ 
\end{theorem}

\vs
\centerline{\Ack Acknowledgments}
\vsp

\noindent
{\ack The authors
would like to thank Prof. Alberto Bressan for suggesting the problem.}

%
%
\newcommand{\auth}{\textsc}
\newcommand{\tit}{\textrm}
\newcommand{\jou}{\textit}

\end{document}